\documentclass[10pt,a4paper,reqno]{amsart}
\usepackage{amsfonts}
\usepackage{amsthm}
\usepackage{amsmath}
\usepackage{mathtools}
\usepackage{dynkin-diagrams}
\usepackage{amscd}
\usepackage[utf8]{inputenc}
\usepackage{t1enc}
\usepackage[mathscr]{eucal}
\usepackage{indentfirst}
\usepackage{graphicx}
\usepackage{graphics}
\usepackage{esint}
\usepackage{pict2e}
\usepackage{epic}
\usepackage{float}
\usepackage{MnSymbol}
\usepackage{multirow}
\usepackage{shuffle} 
\usepackage[table,xcdraw]{xcolor} 
\usepackage{hyperref}
\usepackage{tcolorbox}
\numberwithin{equation}{section}
\usepackage[margin=2.9cm]{geometry}
\usepackage{enumitem}
\usepackage{epstopdf} 
\usepackage{fancyvrb} 
\usepackage{tikz}

\usetikzlibrary{decorations.markings,calc}

\def\centerarc[#1](#2)(#3:#4:#5)
{ \draw[#1] ($(#2)+({#5*cos(#3)},{#5*sin(#3)})$) arc (#3:#4:#5); }

\def\centerarcpath(#1)(#2:#3:#4)
{ ($(#1)+({#4*cos(#2)},{#4*sin(#2)})$) arc (#2:#3:#4) }

\usetikzlibrary{math}

\usepackage[noadjust]{cite}

\newmuskip\pFqmuskip
\newcommand*\pFq[6][8]{%
	\begingroup 
	\pFqmuskip=#1mu\relax
	\mathchardef\normalcomma=\mathcode`,
	\mathcode`\,=\string"8000
	\begingroup\lccode`\~=`\,
	\lowercase{\endgroup\let~}\pFqcomma
	{}_{#2}F_{#3}{\left[\genfrac..{0pt}{}{#4}{#5};#6\right]}%
	\endgroup
}
\newcommand{\pFqcomma}{{\normalcomma}\mskip\pFqmuskip}

\theoremstyle{plain}
\newtheorem{theorem}{Theorem}[section]
\newtheorem{lemma}[theorem]{Lemma}
\newtheorem{corollary}[theorem]{Corollary}
\newtheorem{proposition}[theorem]{Proposition}

\theoremstyle{definition}

\newtheorem{conjecture}[theorem]{Conjecture}
\newtheorem{remark}[theorem]{Remark}

\newtheorem{example}[theorem]{Example}

\newcommand{\ord}[1]{\underset{#1}{\operatorname{\: ord \:}}}

\renewcommand{\Re}{\operatorname{Re}}

\begin{document}
	
	\title{Vanishing of Witten zeta function at negative integers}
	
	\author[Kam Cheong Au]{Kam Cheong Au}
	
	\address{University of Cologne \\ Department of Mathematics and Computer Science \\ Weyertal 86-90, 50931 Cologne, Germany} 
	
	\email{kau@uni-koeln.de}
	

	\begin{abstract} {We introduce a new analytic method for studying Witten zeta function of a root system $\Phi$, based on a refined manipulation of an integral representation involving the Hurwitz zeta function. As an application, we prove high-order vanishing at negative even integers. This technique also describes non-trivially, the arithmetic nature of the leading term, in which the highest root of $\Phi$ makes a surprising appearance.}
	\end{abstract}
	
	\maketitle
	
	\section{Introduction}
	
	Let $G$ be a compact Hausdorff topological group, its \textit{Witten zeta function} (also known as \textit{representation zeta function} \cite{avni2013representation, larsen2008representation}) is defined by $$\zeta_G(s) := \sum_{\rho} \frac{1}{(\text{dim} \: \rho)^s},$$
	where the sum runs over all finite-dimensional irreducible representations of $G$. A familiar example is when $G=\text{SU}(2)$, for which $\zeta_{\text{SU}(2)}(s) = \zeta(s)$, the classical Riemann zeta function. Motivated by simple cases like $G = \text{SU}(3)$ (a Lie group) and $G = \text{SL}_2(\mathbb{Z}_p)$ (a $p$-adic group), Kurokawa and Ochiai \cite{kurokawa2013zeros} conjectured the following:
	\begin{conjecture}
		If $G$ is infinite, then $\zeta_G(-2) = 0$.
	\end{conjecture}
	
	While substantial progress has been made when $G$ is a $p$-adic group~\cite{gonzalez2014representation}, almost nothing is known for Lie groups.
	This article focuses on the latter case and we will show a much stronger statement for compact simply-connected Lie groups. In this setting, such a group 
	$G$ may be identified with its root system $\Phi$, and we write $\zeta_\Phi(s) := \zeta_G(s)$. Our first main result is the following.
	
	\begin{theorem}\label{vanishing_theorem}
		(a) The order of vanishing of $\zeta_\Phi(s)$ at negative even integer is at least the rank of $\Phi$, that is,
		$$\ord{s \in -2\mathbb{N}}(\zeta_\Phi(s)) \geq \text{ Rank of $\Phi$}.$$
		(b) If $\Phi$ is one of $$A_n,\quad D_{n} \text{ with }n\geq 5\text{ and odd},\quad E_6,$$ then $\zeta_\Phi(s)$ vanishes at negative odd integers. More precisely, $$\ord{s \in 1-2\mathbb{N}}(\zeta_\Phi(s)) \geq  \begin{cases}
			\lfloor \frac{n}{2} \rfloor  \quad &\Phi = A_n \text{ for any }n, \\
			1 \quad &\Phi = D_n \text{ with }n\geq 5\text{ and odd},\\
			2 \quad &\Phi = E_6,\\
			0 \quad &\text{otherwise}.\\
		\end{cases}$$
	\end{theorem}
	
	Moreover, we expect that the lower bounds in the theorem above are optimal:
	\begin{conjecture}\label{exact_order_conjecture}
		In Theorem~\ref{vanishing_theorem}, all stated inequalities for orders of vanishing are actually equalities.
	\end{conjecture} 
	
	Witten and Zagier introduced the functions $\zeta_\Phi(s)$ to illustrate their special values at positive even integers~\cite{witten1991quantum, zagier1994values}. A substantial literature has since developed concerning their multivariable generalizations~\cite{komori2023theory, matsumoto2006witten, komori2010witten, komori2015witten, komori2012witten, komori2020zeta, nakamura2006functional}. The standard methods in this area rely on expressing $\zeta_\Phi(s)$ as certain multiple zeta series, for example
	\begin{align}\label{aux_10}
		\zeta_{A_2}(s) &= 2^s \sum_{n,m\geq 1} \frac{1}{m^sn^s(m+n)^s}, \nonumber \\
		\zeta_{B_2}(s) &= 6^s \sum_{n,m\geq 1} \frac{1}{m^sn^s(m+n)^s (m+2n)^s}, \\
		\zeta_{A_3}(s) &= 12^s \sum_{m_1,m_2,m_3\geq 1} \frac{1}{m_1^sm_2^sm_3^s (m_1+m_2)^s(m_2+m_3)^s(m_1+m_2+m_3)^s}. \nonumber
	\end{align}
	These Dirichlet-type representations are closely related to Shintani's work in algebraic number theory~\cite{hida1993elementary}, and many authors have since studied zeta functions of a more general flavor~\cite{friedman2004shintani, bruna2025polynomials, rutard2023values, onodera2021multiple, essouabri2025values, essouabri1997singularite}. \par
	
	Our Theorem \ref{vanishing_theorem} however, would be inaccessible from this conventional view point, largely because the Weyl group symmetry is invisible in this form. Instead, we employ an integral representation (Proposition~\ref{int_rep_of_zeta}), whose integrand involves the Hurwitz zeta function $\zeta(s,a) := \sum_{n\geq 0} (n+a)^{-s}$ and whose prefactors encode the Poincaré polynomial of the root system. For example, \begin{align}\label{aux_11}
		\zeta_{A_2}(s) \: 2^{-s}  \frac{\sin(\pi s)\sin(\frac{3\pi s}{2})}{\sin^2(\frac{\pi s}{2})} &= \left(\frac{(2\pi)^s}{\Gamma(s)}\right)^3 \int_0^1 \zeta(1-s,x) \zeta(1-s,1-x)^2 dx, \nonumber \\
		\zeta_{B_2}(s) \: 6^{-s} \frac{\sin(\pi s)\sin(2\pi s)}{\sin^2(\frac{\pi s}{2})} &= \left(\frac{(2\pi)^s}{\Gamma(s)}\right)^4  \int_{[0,1]^2} \zeta(1-s,x_1)\zeta(1-s,x_2) \\ & \quad\quad\quad\quad\quad\quad\quad\times  \zeta(1-s,\{-x_1-2x_2\})\zeta(1-s,\{-x_1-x_2\}) dx_1 dx_2, \nonumber \end{align}
	with $\{x\}$ the fractional part of $x\in \mathbb{R}$. This representation, while known in the literature~\cite{komori2012witten, komori2023theory}, has not been usefully exploited; in particular, it leads to an especially clean proof of Theorem~\ref{vanishing_theorem}, whose key step is an innovative application of the following purely analytical lemma (Lemma \ref{KeyLemma}). 
	
	\begin{lemma}
		Let $V \subset \mathbb{R}^n$ be an $n$-dimensional simplex, abbreviate $\boldsymbol{x} = (x_1,\cdots,x_n)$. Let $f_1(\boldsymbol{x}),\cdots,f_k(\boldsymbol{x})$ be linear polynomials with real coefficients that are non-zero on interior of $V$; $E_s(\boldsymbol{x})$ be an analytic function on a neighbourhood of $(s,\boldsymbol{x})\in \mathbb{C}\times V$. Then the integral
		\begin{equation}\label{int_Es}\int_V (f_1(\boldsymbol{x})\cdots f_k(\boldsymbol{x}))^{s} E_s(\boldsymbol{x}) d\boldsymbol{x},\qquad \Re(s)>0,\end{equation}
		admits meromorphic extension to $s\in \mathbb{C}$ and the order of any poles is at most $n$.
	\end{lemma}
	
	The two incarnations of Witten zeta function as multiple zeta function (\ref{aux_10}) and definite integral of $\zeta(s,a)$ (\ref{aux_11}), endows it with a very rich theory that have only begun to emerge in recent works~\cite{onodera2014functional, bridges2024asymptotic, romik2017number}. Theorem~\ref{vanishing_theorem} may be just the beginning of a richer landscape, even in rank $2$, it already yields deep conclusions. For example, in \cite[Section~4]{au2024single}, it was shown that vanishing of $\zeta_{A_2}(-2n), \zeta_{B_2}(-2n)$ and $\zeta_{G_2}(-2n)$ are respectively equivalent to following three deep combinatorial identities:
	$$\frac{(2n)!}{(4n+1)!}\zeta(-6n-1) = \sum_{k=0}^{2n} \frac{1}{k! (2 n-k)!} \zeta (-k-2 n) \zeta (k-4 n),$$
	$$ \frac{\zeta(-8n-1)}{3}(1+2^{-1-4n})[f_B(x)^{2n}\log f_B(x)][x^{1+6n}] = \sum_{k=0}^{4n} [f_B(x)^{2n}][x^k] \zeta (-k-2 n) \zeta (k-6 n),$$
	$$\frac{\zeta(-12n-1)}{5}(1+3^{-1-6n})[f_G(x)^{2n}\log f_G(x)][x^{1+10n}] = \sum_{k=0}^{8n} [f_G(x)^{2n}][x^k] \zeta (-k-2 n) \zeta (k-10 n),$$
	here $f_B(x) = (1+x)(1+2x), f_G(x) = (1+x)(1+2x)(1+3x)(2+3x)$ and $[f(x)][x^n]$ means the $x^n$-coefficient in the series expansion of $f(x)$ at $x=0$. Such identities (known as lacunary recurrence) for Bernoulli numbers\footnote{recall Riemann zeta at negative integers are equivalent to Bernoulli numbers} are rare and highly non-trivial \cite{agoh2007convolution}. Among the three identities above, a direct proof is known only for the first one \cite{romik2017number, mertens2015lacunary}.  \par
	After the appearance of this manuscript as a preprint, K.~Onodera communicated to the author that he had independently obtained a result similar to Theorem~\ref{vanishing_theorem} in unpublished work.
	\\[0.02in]
	
	Motivated by the philosophy that the first nonvanishing Taylor coefficient of an $L$-function often carries deep arithmetic information (as in the cases of Dedekind zeta functions or the BSD conjecture), we examine the nature of the quantity
	$$\zeta_\Phi^{(n)}(-m), \qquad m\in 2\mathbb{N}, n = \text{rank}(\Phi),$$
	which, by Conjecture~\ref{exact_order_conjecture}, should be the first nonvanishing coefficient. Although $\zeta_\Phi(s)$ is not an $L$-function\footnote{for example, $\zeta_\Phi(s)$ lacks a functional equation and an Euler product}, this number still has nice arithmetic structure, with a surprising connection to the highest root of $\Phi$. This is our second main result. 
	
	For an irreducible $\Phi$, denote $\mathcal{H}(\Phi) \subset \mathbb{N}$ to be the set of coefficients appearing when expressing its highest root in terms of simple roots. 
	\begin{theorem}\label{special_value_thm}
		For irreducible root system $\Phi$ of rank $n$, let $r$ be the number of positive roots in $\Phi$, let $T$ be the set of rational numbers in $(0,1]$ whose denominator is an element of $\mathcal{H}(\Phi)$. Let $m\in 2\mathbb{N}$, then $$\pi^{mr} \zeta_\Phi^{(n)}(-m)$$
		is in the $\mathbb{Q}$-span of numbers
		\begin{equation}\tag{*}\label{aux_6}\left\{\zeta(m_1,t_1)\zeta(m_2,t_2)\cdots\zeta(m_k,t_k) \bigg| 1\leq k\leq n,\quad m_j\in \mathbb{N},\quad \sum_{j=1}^k m_j = mr,\quad  m_j\geq m+1, \quad t_j\in T \right\}.\end{equation}
	\end{theorem}
	
	In particular, for classical root systems, the second argument of the Hurwitz zeta function is not needed (Remark after Corollary~\ref{coeff_form_corollary}), yielding the following corollary.
	\begin{corollary}
		When $\Phi$ is one of $A_n, B_n, C_n$ or $D_n$, let $r$ be the number of positive roots in $\Phi$, $m\in 2\mathbb{N}$, then $\pi^{mr} \zeta_\Phi^{(n)}(-m)$ is a polynomial in $\pi^2$ and $\{\zeta(2j+1)\}_{j\in \mathbb{N}}$. 
	\end{corollary}
	
	There also exists corresponding statement for odd $m$ (Corollary \ref{coeff_form_corollary}). For both Theorems~\ref{vanishing_theorem} and \ref{special_value_thm}, the sole case previously known in the literature is $\Phi = A_2$. In this case, Onodera \cite{onodera2014functional, onodera2021multiple} derived an explicit formula: when $m\in 2\mathbb{N}$, \begin{multline}\label{A2leadingcoeff}\frac{(2\pi)^{3m}}{(m!)^3} 2^{m} \zeta''_{A_2}(-m) = \frac{(-1)^{m/2}}{3} \bigg[\sum_{j=0}^m \binom{m+j}{j}\binom{2m-j}{m-j} \zeta(1+m+j)\zeta(1+2m-j) (1+2(-1)^j) \\ + \frac{1}{2}\binom{3m+1}{2m+1} \zeta(3m+2)\bigg].\end{multline}
	We provide a more conceptual derivation of his identity (Remark~\ref{A2_alternative_proof}), illustrating the method used to prove Theorem~\ref{special_value_thm} in the simple case $\Phi = A_2$. The proof in the general case combines analytic and combinatorial ingredients. The analytic component consists of a refinement of the previous lemma, while the combinatorial component exploits symmetry of the root system and it is where $\mathcal{H}(\Phi)$ enters naturally into the picture. \par
	
	More precisely, the appearance of $\mathcal{H}(\Phi)$ arises from the following root-system-theoretic result (Proposition~\ref{EH_prop}) that relates two quantities not obviously seen to be connected, and strengthens earlier work of Springer and Steinberg \cite{Springer1970} on bad primes of $\Phi$ (Remark \ref{springer_remark}).
	
	\begin{proposition}
		Let $\Phi$ be an irreducible root system in Euclidean space $V$, let $L(\Phi)$ be its root lattice. Define\footnote{recall exponent of a group is the smallest positive integer that annihilates the group.} $$\mathcal{E}(\Phi) := \left\{\text{Exponent of the finite abelian group }\frac{L(\Phi)} {\text{Span}_\mathbb{Z} S} \biggr| S\subset \Phi \text{ such that }\text{Span}_\mathbb{R}S=V\right\} \subset \mathbb{N}.$$
		Then $\mathcal{E}(\Phi) = \mathcal{H}(\Phi) \cup \{1\}$.
	\end{proposition}

	We structure the article as follows. In Section~2 we derive the integral representation of $\zeta_\Phi(s)$ involving the Hurwitz zeta function, which forms the backbone of the subsequent arguments.
	Section~3 contains the proof of Theorem~\ref{vanishing_theorem}, assuming an analytic lemma proved in Section~4.
	Section~5 states and gives examples for Theorem~\ref{special_value_thm}, whose proof occupies Sections~6 and~7: the former analytic in nature, the latter combinatorial in nature.

	\section*{Acknowledgment}
	The author thanks Prof. Kazuhiro Onodera for valuable feedback. The author has received funding from the European Research Council (ERC) under the European Union’s Horizon 2020 research and innovation programme (grant agreement No. 101001179).

	\section{An integral representation of $\zeta_\Phi(s)$}
	We start by deriving an integral representation of $\zeta_\Phi(s)$ that is essential for us. Throughout this article, we adopt following notations attached to a root system $\Phi$:
	\begin{itemize}
		\item $n := \text{rank}(\Phi)$, with inner product pairing $(\cdot,\cdot)$ and coroot $\alpha^\vee := \frac{2\alpha}{(\alpha,\alpha)}$; 
		\item $\Lambda$: weight lattice associated to $\Phi$;
		\item $\Delta := \{\alpha_1,\cdots,\alpha_n\}$: a fixed set of simple roots;
		\item $\Phi^+$: positive roots with respect to $\Delta$; 
		\item $r := |\Phi^+|$, the number of positive roots; 
		\item $\Lambda^+$: dominant weights in $\Lambda$ with respect to $\Delta$;
		\item $W$: Weyl group of $\Phi$;
		\item $l(\sigma)$ with $\sigma \in W$: length of a Weyl group element.
	\end{itemize}
	
	Let $\lambda_1,\cdots,\lambda_n$ to be fundamental weights associated to $\Delta$, they are characterized by
	$$(\lambda_i,\alpha_j^\vee) = \begin{cases}1 \quad &i=j\\ 0 \quad &i\neq j\end{cases}.$$
	Denote also $$\delta := \frac{1}{2}\sum_{\alpha\in \Phi^+} \alpha = \sum_{i=1}^n \lambda_i,$$
	this is the smallest strongly dominant weight, the set of strongly dominant weights equals $\delta + \Lambda^+$. 
	
	Finite-dimensional irreducible representation of $\mathfrak{g}$ is in one-one correspondence to dominant weight $\lambda \in \Lambda^+$. For irreducible representation associated with $\lambda \in \Lambda^+$, Weyl dimension formula implies its dimension is
	$$\frac{\prod_{\alpha\in \Phi^+}(\lambda+\delta,\alpha^\vee)}{\prod_{\alpha\in \Phi^+}(\delta,\alpha^\vee)}.$$
	We shall denote $K_\Phi := \prod_{\alpha\in \Phi^+}(\delta,\alpha^\vee) \in \mathbb{N}$, whose value is not important to us. We can write $\zeta_\Phi(s)$ as
	$$\zeta_\Phi(s) = \sum_{\pi} \frac{1}{(\dim \pi)^s} = K_\Phi^s \sum_{\lambda \text{ strongly dominant}} \frac{1}{\prod_{\alpha\in \Phi^+} (\lambda,\alpha^\vee)^s}.$$
	We are going to complete the above sum to all $\lambda\in \Lambda$. First we agree, for a negative real number $r$, \begin{equation}\label{negative_base}r^s = |r|^s e^{-\pi r s} \qquad r<0, s\in \mathbb{C}.\end{equation} Consider the completed sum: \begin{equation}\label{completed_sum}{\sum_{\lambda \in \Lambda}}' \frac{1}{\prod_{\alpha\in \Phi^+} (\lambda,\alpha^\vee)^s},\end{equation}
	here we are summing over all $\lambda \in \Lambda$, but excluding those $\lambda$ such that the denominator vanishes. Geometrically, we are summing $\lambda$ in the interior of each Weyl chamber, on which the Weyl group $W$ acts freely and transitively. Therefore (\ref{completed_sum}) equals
	$$\sum_{\sigma\in W} \sum_{\lambda \text{ strongly dominant}} \frac{1}{\prod_{\alpha\in \Phi^+} (\sigma^{-1}\lambda,\alpha^\vee)^s} = \sum_{\sigma\in W} \sum_{\lambda \text{ strongly dominant}} \frac{1}{\prod_{\alpha\in \Phi^+} (\lambda,(\sigma\alpha)^\vee)^s}.$$
	For each $\alpha \in \Phi^+$, $\sigma\alpha$ could be either positive or negative; number of positive roots that are sent to negative roots by $\sigma$ is exactly $l(\sigma)$, so $$\prod_{\alpha\in \Phi^+} (\lambda,(\sigma\alpha)^\vee)^s = e^{-\pi i l(\sigma) s}\prod_{\alpha\in \Phi^+} (\lambda,\alpha^\vee)^s.$$
	Writing $$P_\Phi(s) := \sum_{\sigma\in W} e^{\pi i l(\sigma) s}, \qquad \xi_\Phi(s) := K_\Phi^{-s} \zeta_\Phi(s),$$
	the exponential polynomial $P_\Phi(s)$ is known as \textit{Poincaré polynomial} of $\Phi$. Expression (\ref{completed_sum}) then equals $P_\Phi(s) \xi_\Phi(s).$
	
	Next we derive an integral representation of (\ref{completed_sum}). Introduce function $$F(s,a) := \sum_{n\neq 0} \frac{e^{2\pi i n a}}{n^s} = \sum_{n\geq 1} \frac{e^{2\pi i n a}}{n^s} + e^{\pi i s} \sum_{n\geq 1}\frac{e^{-2\pi i n a}}{n^s},\qquad \Re(s)>1, a\in \mathbb{R},$$
	here we remind readers of our convention (\ref{negative_base}). Note that $$n^{-s} = \int_0^1 e^{-2\pi i n x} F(s,x) dx, \qquad n\neq 0, n\in \mathbb{Z}.$$ 
	
	When we write $\lambda = m_1\lambda_1+\cdots+m_n \lambda_n$ as linear combination of fundamental weights, $(\lambda,\alpha_j^\vee) = m_i$, expression (\ref{completed_sum}) becomes
	$$\begin{aligned}
		&\quad \sum_{m_i\neq 0} m_1^{-s}\cdots m_n^{-s} \prod_{\alpha \in \Phi^+ - \Delta} \left(\sum_{i=1}^n m_i (\lambda_i,\alpha^\vee)\right)^{-s} \\
		&= \sum_{m_i\neq 0} m_1^{-s}\cdots m_n^{-s} \prod_{\alpha \in \Phi^+ - \Delta} \left(\int_0^1 \exp\left[-2\pi i x_\alpha \sum_{i=1}^n m_i (\lambda_i,\alpha^\vee)\right] F(s,x_\alpha) dx_\alpha \right) \\
		&= \int_{[0,1]^{i-n}}\prod_{i=1}^n F\left(s,-\sum_{\alpha \in \Phi^+ - \Delta} (\lambda_i, \alpha^\vee) x_\alpha\right) \prod_{\alpha \in \Phi^+ - \Delta} F(s,x_\alpha) dx_\alpha,
	\end{aligned}$$
	here we assume $\Re(s)$ sufficiently large (actually $\Re(s)>1$ is enough) such that interchanges of sum and integral in the last step are valid. Note the final integral has dimension $|\Phi^+| - |\Delta| = r - n$. \par It is not obvious from the definition of $F(s,a)$ that it can be extended to a complex analytic function in $a$, but this is indeed true:
	\begin{lemma}\cite[p.~257]{apostol2013introduction}
		For $\Re(s)>1$ and $0<a<1$, 
		$$\zeta(1-s,a) = \frac{\Gamma(s)}{(2\pi i)^s} F(s,a).$$
	\end{lemma}
	
	Combining the two expressions of equation (\ref{completed_sum}), we obtain
	\begin{proposition}\label{int_rep_of_zeta}
		For $\Re(s)>1$, let \begin{equation}\label{int_rep_I}I_\Phi(s) :=  \int_{[0,1]^{r-n}} \prod_{\alpha \in \Phi^+ - \Delta} \zeta(1-s,x_\alpha) \prod_{i=1}^n \zeta\left(1-s,\left\{ -\sum_{\alpha \in \Phi^+ - \Delta} (\lambda_i, \alpha^\vee) x_\alpha \right\}\right) dx_\alpha,\end{equation} where $\{x\}$ means fractional part of a real number $x$. Then \begin{equation}\label{expr_for_zetaPhi}P_\Phi(s) \xi_\Phi(s) = \left(\frac{(2\pi i)^s}{\Gamma(s)}\right)^i I_\Phi(s).\end{equation}
	\end{proposition}
	We remark that $\Re(s)>1$ is not the largest region in which equation (\ref{int_rep_I}) holds, but this is not important to us. 
	
	As examples, here are expressions of $I_\Phi(s)$ for some low rank root systems:
	$$\begin{aligned}
		I_{A_2}(s) &= \int_0^1 \zeta(1-s,x) \zeta(1-s,1-x)^2 dx, \\
		I_{B_2}(s) &= \int_{[0,1]^2} \zeta(1-s,x_1)\zeta(1-s,x_2) \zeta(1-s,\{-x_1-2x_2\})\zeta(1-s,\{-x_1-x_2\}) dx_i,  \\
		I_{A_3}(s) &= \int_{[0,1]^3} \prod_{i=1}^3 \zeta(1-s,x_i) \zeta(1-s,\{-x_1-x_3\})\zeta(1-s,\{-x_2-x_3\})\zeta(1-s,\{-x_1-x_2-x_3\}) dx_i, \\ 
		I_{G_2}(s) &= \int_{[0,1]^4} \prod_{i=1}^4 \zeta(1-s,x_i) \zeta(1-s,\{-x_1-x_2-2x_3-x_4\})\zeta(1-s,\{-3x_1-x_2-3x_3-2x_4\}) dx_i.
	\end{aligned}$$
	
	\begin{remark}
		Formula (\ref{expr_for_zetaPhi}) itself is nothing new \cite{komori2012witten},\cite[Chapter~4]{komori2023theory}, it provides a quick proof that $\zeta_\Phi(2m) \in \mathbb{Q} \pi^{2mr}$ for positive integer $m$: using $\zeta(-n,a) = -\frac{B_{n+1}(a)}{n+1}$ with $B_n(x)$ the Bernoulli polynomial, we see $I_\Phi(2n) \in \mathbb{Q}$ because integrand is a $\mathbb{Q}$-polynomial and region of integration is a finite union of  $\mathbb{Q}$-simplex, and $P_\Phi(2n) = \sum_{\sigma\in W} e^{2 \pi i l(\sigma) n} = |W| \neq 0$. Zagier's original proof \cite{zagier1994values}, based on partial fraction, is somewhat different.
	\end{remark}
	
	\begin{remark}
		The Poincaré polynomial admits a nice factorization \cite[Chapter~3]{humphreys1992reflection} \begin{equation}\label{aux_3}P_\Phi(s) = \sum_{\sigma\in W}e^{i\pi l(\sigma)s} = \prod_{k=1}^n \frac{e^{\pi i s d_k}-1}{e^{\pi i s}-1},\end{equation}
		for some positive integers $d_1,\cdots,d_n$, known as \textit{degrees} of $W$. We tabulate their values below. Using this factorization, we obtain an alternative form of equation (\ref{expr_for_zetaPhi}):
		\begin{equation}\label{sine_form}\xi_\Phi(s) \prod_{k=1}^n \frac{\sin(\frac{\pi s d_k}{2})}{\sin(\frac{\pi s}{2})} = \left( \frac{(2\pi)^s}{\Gamma(s)}\right)^r I_\Phi(s).\end{equation}
		\vspace*{-3mm}
		\begin{table}[h]
			\centering
			\begin{tabular}{c|c}
				$\Phi$ & $d_1,\cdots,d_n$ \\ \hline
				$A_n$ & $2,3,\cdots,n+1$ \\
				$B_n, C_n$ & $2,4,6,\cdots,2n$ \\
				$D_n$ & $2,4,6,\cdots,2n-2,n$ \\
				$G_2$ & $2,6$ \\
				$F_4$ & $2,6,8,12$ \\
				$E_6$ & $2,5,6,8,9,12$ \\
				$E_7$ & $2,6,8,10,12,14,18$ \\
				$E_8$ & $2,8,12,14,18,20,24,30$
			\end{tabular}
			\caption{\small Degree of Weyl group associated with irreducible root systems.}
			\label{Weyl_degree:degree}
		\end{table}
	\end{remark}
	\vspace*{-8mm}
	
	\section{A lemma and deduction of Theorem \ref{vanishing_theorem}}
	We shall prove the following lemma in the next section. \begin{lemma}\label{KeyLemma}
		Let $V \subset \mathbb{R}^n$ be an $n$-dimensional simplex, abbreviate $\boldsymbol{x} = (x_1,\cdots,x_n)$. Let $f_1(\boldsymbol{x}),\cdots,f_k(\boldsymbol{x})$ be linear polynomials over $\mathbb{R}$ that are non-zero on interior of $V$; $E_s(\boldsymbol{x})$ be an analytic function on a neighbourhood of $(s,\boldsymbol{x})\in \mathbb{C}\times V\subset \mathbb{C}^{n+1}$. Then the integral
		\begin{equation}\label{int_Es}\int_V (f_1(\boldsymbol{x})\cdots f_k(\boldsymbol{x}))^{s} E_s(\boldsymbol{x}) dx_1\cdots dx_n,\qquad \Re(s)>0,\end{equation}
		\begin{enumerate}[leftmargin=*]
			\item admits meromorphic extension to $\mathbb{C}$;
			\item all poles are at $s\in \mathbb{Q}$;
			\item order of poles at any point is $\leq n$.
		\end{enumerate}
	\end{lemma}
	
	\begin{remark}
	The analytic continuation of such integrals has been studied in the context of algebraic $D$-module theory and Bernstein–Sato polynomials~\cite{bernshtein1968feasibility, kashiwara1976b, coutinho1995primer, kashiwara2003d}, but we provide a  self-contained argument. This is partly because the lemma is by no means obvious\footnote{points (1), (2) are known, but (3) is more non-trivial} even with this general theory in hand, and partly because we later require an algebraic version (Lemmas~\ref{V_lemma} and \ref{V_defined_lemma}) to establish our second main result (Theorem \ref{special_value_thm}), for which $D$-module methods no longer apply.
	\end{remark}
	
	\begin{corollary}
		For a root system $\Phi$, $I_\Phi(s)$ has a meromorphic continuation to $\mathbb{C}$, order of pole at any point $s\neq 0$ is $\leq r-n$. 
	\end{corollary}
	\begin{proof}
		In the expression
		$$I_\Phi(s) =  \int_{[0,1]^{r-n}} \prod_{\alpha \in \Phi^+ - \Delta} \zeta(1-s,x_\alpha) \prod_{i=1}^n \zeta\left(1-s,\left\{ -\sum_{\alpha \in \Phi^+ - \Delta} (\lambda_i, \alpha^\vee) x_\alpha \right\}\right) dx_\alpha$$
		the fractional part $\{\cdots\}$ is piecewise continuous, we can triangulate  integration domain $[0,1]^{r-n}$ into finitely smaller $(r-n)$-simplex $V_i$, $I_\Phi(s) = \sum_i I_{V_i}(s)$, where $I_{V_i}(s)$ is the integration restricted on $V_i$ such that on each of $V_i$, the second argument of Hurwitz zeta functions are linear polynomials whose values on $V_i$ are between $0$ and $1$, but may attain $0$ on its boundary. Fix an $i$, abbreviate $V_i$ as $V$, write
		$$I_{V}(s) = \int_{V} \zeta(1-s,f_1(\boldsymbol{x})) \cdots \zeta(1-s,f_r(\boldsymbol{x})) d\boldsymbol{x}$$
		with $f_i$ being linear polynomials. Note that $s\: \zeta(s,a)$ is analytic on $(s,a)\in \mathbb{C}\times \{a : \Re(a)>0\}$, with $a=0$ being a singularity. In the above integral, the place where $f_r(\boldsymbol{x})$ vanishes, which could happen on the boundary of $V$, is problematic. To tame this issue, use the identity $$\zeta(1-s,x) = \zeta(1-s,1+x) + x^{s-1}$$ and rewrite $I_V(s)$ as
		$$s^r I_{V}(s) = \int_{V} [s\zeta(1-s,1+f_1(\boldsymbol{x})) + s f_1(\boldsymbol{x})^{s-1}] \cdots [s\zeta(1-s,1+f_r(\boldsymbol{x})) + s f_r(\boldsymbol{x})^{s-1}] d\boldsymbol{x}.$$
		Now $1+f_i(\boldsymbol{x})$ is positive on a neighbourhood of $V$, $s \zeta(1-s,1+f_i(\boldsymbol{x}))$ is thus analytic functions on a neighbourhood of $(s,\boldsymbol{x}) \in \mathbb{C}\times V$. The integrand can be expanded into $2^r$ terms, each according to Lemma \ref{KeyLemma}, admits meromorphic continuation to $\mathbb{C}$, and order of poles at $s\neq 0$ is at most dimension of integral, i.e. $r-n$. Therefore their sum $s^r I_{V}(s)$ has this property too, so is $s^r I_\Phi(s)$.
	\end{proof}
	
	\begin{proof}[Proof of Theorem \ref{vanishing_theorem}]
		Recall formula (\ref{expr_for_zetaPhi}), $$P_\Phi(s)  \xi_\Phi(s) = \left(\frac{(2\pi i)^s}{\Gamma(s)}\right)^r I_\Phi(s).$$
		When $s$ is a negative even integer, $P_\Phi(s) = \sum_{\sigma\in W} e^{\pi i l(\sigma) s} = |W| \neq 0$, the gamma factor $\left(\frac{(2\pi i)^s}{\Gamma(s)}\right)^r$ has zero of order $r$, and $I_\Phi(s)$ can have poles of order at most $r-n$, thus $\zeta_\Phi(s)$ has zero of order at least $n$.\\
		When $s$ is a negative odd integer, the factorization (\ref{aux_3}) implies
		$$\ord{s\in 2\mathbb{Z}+1} P_\Phi(s) = \text{number of even degrees}.$$
		Consequently, \begin{align*}\ord{s\in -2\mathbb{N}+1}\zeta_\Phi(s) &= -\left(\ord{s\in -2\mathbb{N}+1}P_\Phi(s) \right)+ \left(\ord{s\in -2\mathbb{N}+1}I_\Phi(s) \right) + r \\ &\geq  -\left(\ord{s\in -2\mathbb{N}+1}P_\Phi(s) \right) + n \\ &= \text{number of odd degrees}.\end{align*}
		The statement then follows by inspecting Table \ref{Weyl_degree:degree}.
	\end{proof}

	\section{Proof of Lemma \ref{KeyLemma}}
	Now we prove Lemma \ref{KeyLemma}, we do it by induction on $n$, i.e. dimension of the integral. Note the space of functions satisfying the three properties in the statement is closed under addition and scalar multiplication. \par
	\textit{Step 1: the base case }$n=1$. We can assume $V=[0,1]$, $f(x)$ is a polynomial with no roots in its interior $(0,1)$. We can assume $f(x) = x^A (1-x)^B$, since any other factor, when raised to $s$-th power, is analytic in a neighbourhood of $V$ and hence can be absorbed into the term $E_s(x)$. Therefore we need to prove
	$$\int_0^1 x^{As} (1-x)^{Bx} E_s(x) dx,\quad A,B\in \mathbb{N}$$
	satisfies the three conditions. Split the integration interval into $[0,1/2]$ and $[1/2,1]$, make $x\mapsto x/2$ in former, $x\mapsto 1-x/2$ in latter, above is
	$$\frac{1}{2}\int_0^1 (\frac{x}{2})^{As} (1-\frac{x}{2})^{Bs} E_s(\frac{x}{2}) dx + \frac{1}{2} \int_0^1 (1-\frac{x}{2})^{As} (\frac{x}{2})^{Bs} E_s(1-\frac{x}{2}) dx.$$
	In former integral, we can absorb $(1-x/2)^{Bs}$ into $E_s(x/2)$; similar for latter integral. Hence it suffices to prove $$\int_0^1 x^{As} E_s(x) dx, \quad A\in \mathbb{N}$$ satisfies the three conditions, where $E_s(x)$ is analytic on a neighbourhood of $(s,x) \in \mathbb{C}\times [0,1]$. If $A = 0$, the integral is entire in $s$, so we assume $A>0$. By agreeing $0\leq \arg x < 2\pi$, $C$ be the contour in figure, we have $$\int_0^1 x^{As} E_s(x) dx  = (e^{2\pi i A s}-1)^{-1} \int_C x^{As} E_s(x) dx, \qquad \Re(s)>0.$$ 		\begin{figure}[h]
		\centering
		\begin{tikzpicture}[decoration={markings,
				mark=at position 0.3cm with {\arrow[line width=1pt]{>}},
				mark=at position 5cm with {\arrow[line width=1pt]{>}},
				mark=at position 7.85cm with {\arrow[line width=1pt]{>}},
				mark=at position 13cm with {\arrow[line width=1pt]{>}}
			}
			]
			\draw[help lines,->] (-1,0) -- (8,0) coordinate (xaxis);
			\draw[help lines,->] (0,-1) -- (0,1) coordinate (yaxis);

			\path[draw,line width=0.8pt,postaction=decorate] (7.5,0) -- (7.5,0.3) -- (0.0,0.3) \centerarcpath(0,0)(90:270:0.3) -- (7.5,-0.3) -- (7.5,0);
			
			\filldraw[black] (7.5,0) circle (1pt) node[anchor=west]{$1$};
		\end{tikzpicture}\caption{\small \label{figureC} The contour $C$.}
	\end{figure}
	
	As $C$ does not pass through $x=0$, $\int_C x^{As} E_s(x) dx$ is an entire function in $s$, $(1-e^{2\pi i A s})^{-1}$ has simple poles at $s\in \frac{1}{A}\mathbb{Z}$, this shows $(*)$ is true when $n=1$. The base case of induction is established. \par
	\textit{Step 2: Barycentric subdivision of a simplex}. For general $n$, we need a suitable triangulation for our $n$-dimensional simplex $V$.
	
	\begin{lemma}
		Let $V$ an $n$-dimensional simplex in $\mathbb{R}^n$, $\mathcal{H}$ be a family of hyperplanes in $\mathbb{R}^n$ such that each its member is disjoint from interior of $V$. Then we can triangulate $V$ into finitely many $n$-simplexes $V = \bigcup V_i$, each $V_i$ has a distinguished vertex $p_i$ such that for any $i$, each hyperplane in $\mathcal{H}$ either contains $p_i$ or is disjoint from $V_i$.
	\end{lemma}
	\begin{proof}
		We claim the barycentric subdivision \cite[p.~119-120]{hatcher2002algebraic} of $V$ satisfies the requirement. For $0\leq j \leq n$, let $\mathcal{F}_i$ be set of $j$-dimensional faces of $V$, for example $\mathcal{F}_0$ consists of vertices in $V$, $\mathcal{F}_1$ consists of lines of $V$, $\mathcal{F}_n = \{V\}$, write $\mathcal{F} = \cup_{j=0}^n \mathcal{F}_j$, note that $|\mathcal{F}| = 2^{n+1}-1$. It is easy to show barycentric subdivision $V = \bigcup V_i$ satisfies following properties
		\begin{itemize}
			\item Each $V_i$ contains a unique element of $\mathcal{F}_0$, say $p_i$,
			\item Each element of $\mathcal{F}$ either contains $p_i$ or is disjoint from $V_i$. 
		\end{itemize}
		
		\begin{figure}[h]
			\centering
			\scalebox{1}{
				\begin{tikzpicture}
					\tikzmath{\x1 = -4; \y1 =-0.9;
						\x2 = 1.6; \y2 = -1.2;
						\x3 = 0; \y3 = 4; 
						\x4 = (\x1 + \x2)/2; \y4 = (\y1 + \y2)/2;
						\x5= (\x1 + \x3)/2; \y5= (\y1 + \y3)/2;
						\x6 = (\x2 + \x3)/2; \y6 = (\y2 + \y3)/2;}
					\draw (\x1,\y1)--(\x2,\y2);
					\draw (\x2,\y2)--(\x3,\y3);
					\draw (\x1,\y1)--(\x3,\y3);
					\draw (\x1,\y1)--(\x6,\y6);
					\draw (\x2,\y2)--(\x5,\y5);
					\draw (\x3,\y3)--(\x4,\y4);
					\filldraw[black] (\x1,\y1) circle (1pt) node[anchor=north]{\small $p_1=p_2$};
					\filldraw[black] (\x2,\y2) circle (1pt) node[anchor=north]{\small $p_3=p_4$};
					\filldraw[black] (\x3,\y3) circle (1pt) node[anchor=south]{\small $p_5 = p_6$};
					\filldraw[black] (-2.0143,0.242536) node[anchor=south]{$V_1$};
					\filldraw[black] (0.6143,0.142536) node[anchor=south]{$V_4$};
					\filldraw[black] (-0.143,-0.642536) node[anchor=south]{$V_3$};
					\filldraw[black] (-1.743,-0.642536) node[anchor=south]{$V_2$};
					\filldraw[black] (-1.043,1.7) node[anchor=south]{$V_6$};
					\filldraw[black] (0.143,1.6) node[anchor=south]{$V_5$};
			\end{tikzpicture}}
			\caption{\small Barycentric subdivision of a $2$-simplex, with $V_i$ and $p_i$ labeled.}
		\end{figure}
		By assumption, each $H \in \mathcal{H}$ is disjoint from the interior of $V$, this forces $H\cap V$ to be either empty or an element of $\mathcal{F}$. Thus $H \cap V$ either contains $p_i$ or is disjoint from $V_i$ according to the second bullet, as claimed. 
	\end{proof}
	
	\textit{Step 3: Induction step}. Let $f(\boldsymbol{x}) := f_1(\boldsymbol{x})\cdots f_k(\boldsymbol{x})$, where $f_i$ as in statement of Lemma \ref{KeyLemma}, let $\mathcal{H}$ consists of zero sets of $f_1,\cdots,f_k$. By hypothesis, they are disjoint from interior of $V$. Using the triangulation of $V$ described in step 2, write $$\int_V (f(\boldsymbol{x}))^{s} E_s(\boldsymbol{x}) d\boldsymbol{x}= \sum_i \int_{V_i} (f(\boldsymbol{x}))^{s} E_s(\boldsymbol{x}) d\boldsymbol{x},$$ it suffices to show integral over each $V_i$ satisfies the three conditions. Focusing on one such $V_i$, we can write $$f(\boldsymbol{x}) = f_{i_1}(\boldsymbol{x}) \cdots f_{i_r}(\boldsymbol{x}) q(\boldsymbol{x}),$$
	where linear polynomials $f_{i_1},\cdots,f_{i_r}$ vanish at our distinguished vertex $p_i$; $q$ is non-zero on $V_i$, then $E_{i,s}(\boldsymbol{x}) := E_s(\boldsymbol{x}) q(\boldsymbol{x})^s$ is analytic on a neighbourhood of $(s,\boldsymbol{x})\in \mathbb{C}\times V_i$. Then $$\int_{V_i} (f(\boldsymbol{x}))^{s} E_s(\boldsymbol{x}) d\boldsymbol{x} = \int_{V_i} f_{i_1}(\boldsymbol{x})^s \cdots f_{i_r}(\boldsymbol{x})^s E_{i,s}(\boldsymbol{x}) d\boldsymbol{x}.$$
	Translating the point $p_i$ to origin, then applying a suitable linear transformation will map $V_i$ to the standard $n$-simplex 
	$$S_n := \{(x_1,\cdots,x_n) \in \mathbb{R}^n | x_1\geq 0, \cdots, x_n\geq 0, 0\leq x_1+\cdots+x_n\leq 1\}.$$
	After this transformation, $f_{i_1},\cdots,f_{i_r}$ will be linearly homogeneous (since they vanish at $p_i$). 
	
	\par It remains to show the three conditions hold for integrals of form $$I(s) := \int_{S_n} f_1(\boldsymbol{x})^s \cdots f_k(\boldsymbol{x})^s E_s(\boldsymbol{x}) d\boldsymbol{x},$$
	where $f_i$ are homogeneous linear and $E_s(\boldsymbol{x})$ analytic on a neighbourhood of $(s,\boldsymbol{x})\in \mathbb{C}\times S_n$. To $I(s)$, we perform the substitution $(x_1,\cdots,x_n) \mapsto (u,y_1,\cdots,y_{n-1}),$
	\begin{equation}\label{int_sub}
		x_1 = uy_1, \quad x_2 = uy_2, \quad \cdots \quad x_{n-1} = uy_{n-1}, \quad x_n = u(1-y_1-\cdots-y_{n-1}), 
	\end{equation}
		whose Jacobian is $u^{n-1}$. Write $\boldsymbol{y} = (y_1,\cdots,y_{n-1})$, we have $f_i(\boldsymbol{x}) = u g_i(\boldsymbol{y})$ with $g_i$ a linear polynomial (not necessarily homogeneous) in $\boldsymbol{y}$. Then $$I(s) = \int_{S_{n-1}} g_1(\boldsymbol{y})^s \cdots g_k(\boldsymbol{y})^s  \left(\int_0^1 u^{ks+n-1} E_s(u \boldsymbol{y}, u(1-y_1-\cdots-y_{n-1})) du\right) d\boldsymbol{y}.$$
	Let $$\begin{aligned}F_s(\boldsymbol{y}) &:= \frac{e^{2\pi i k s}-1}{2\pi i}\int_0^1 u^{ks + n-1} E_s( u \boldsymbol{y}, u(1-y_1-\cdots-y_{n-1})) du \\ 
		&= \frac{1}{2\pi i} \int_C u^{ks + n-1} E_s( u \boldsymbol{y}, u(1-y_1-\cdots-y_{n-1})) du,\end{aligned}$$
	where $C$ is the contour used above. Because $C$ does not pass through $0$, we see $F_s(\boldsymbol{y})$ is analytic in a neighbourhood of $(s,\boldsymbol{y})\in \mathbb{C}\times S_{n-1}$. Finally, $$I(s) = \frac{2\pi i}{e^{2\pi i ks}-1} \int_{S_{n-1}} g_1(\boldsymbol{y})^s \cdots g_k(\boldsymbol{y})^s F_s(\boldsymbol{y}) d\boldsymbol{y}.$$
	The integral now has a lower dimension, so we can apply the induction hypothesis on it. The factor $(e^{2\pi i ks}-1)^{-1}$ increases order of pole by $1$, so $I(s)$ has poles of order at most $n$, completing the proof of Lemma \ref{int_Es}.
	
	\begin{corollary}\label{orderofpolelessthank}
		Using notations in the statement of Lemma \ref{KeyLemma}, the order of poles of the integral is also $\leq k$, the number of linear forms appearing.
	\end{corollary}
	\begin{proof}
		If $k\geq n$, there is nothing to prove. Assume $k<n$, $\boldsymbol{x} = (x_1,\cdots,x_k,x_{k+1},\cdots,x_n)$, after an invertible linear transformation, we can assume $f_1,\cdots,f_k$ to be independent of $x_{k+1},\cdots,x_n$. Let $\pi: \mathbb{R}^n \to \mathbb{R}^k$ be projection onto first $k$ coordinates. Write $\boldsymbol{z} = (x_1,\cdots,x_k), \boldsymbol{w} = (x_{k+1},\cdots,x_n)$, our integral in question is \begin{equation}\label{aux3}\int_V (f_1(\boldsymbol{z})\cdots f_k(\boldsymbol{z}))^s E_s(\boldsymbol{z},\boldsymbol{w}) d\boldsymbol{z}d\boldsymbol{w}.\end{equation}
		We need to prove its poles are of order $\leq k$. Above expression equals (Fubini's theorem), 
		\begin{equation}\label{aux4}\int_{\pi(V)} (f_1(\boldsymbol{z})\cdots f_k(\boldsymbol{z}))^s G_s(\boldsymbol{z}) d\boldsymbol{z},\end{equation}
		where $$G_s(\boldsymbol{z}) = \int_{V\cap \pi^{-1}(\boldsymbol{z})} E_s(\boldsymbol{z},\boldsymbol{w}) d\boldsymbol{w}.$$
		$G_s(\boldsymbol{z})$ is analytic on a neighbourhood of $(s,\boldsymbol{z}) \in \mathbb{C}\times \pi(V)$ because $E_s(\boldsymbol{x})$ is analytic on a neighbourhood of $\mathbb{C}\times V$. Also $\pi(V)$ is a union of $k$-dimensional simplexes, Lemma (\ref{KeyLemma}) applied to integral (\ref{aux4}) says its order of poles is $\leq k$.
	\end{proof}

	\section{Leading term of $I_\Phi(s)$ at negative integer}
	Theorem \ref{rationality_theorem} below, whose proof will occupy us for the next two sections, provides the arithmetical nature of leading coefficient of $I_\Phi(s)$ at a negative integer. \par
	
	We need first some notations. For a meromorphic function $f(s)$ and $s_0\in \mathbb{C}, k\in \mathbb{Z}$, let $[f(s)][(s-s_0)^{k}]$ be its $(s-s_0)^k$-coefficient of $f(s)$ at $s=s_0$.\par
	For $m,l,w\in \mathbb{N}$, $T$ a set of positive rational numbers, define
	$$\mathcal{Z}(T)_{m,l,w} := \text{Span}_\mathbb{Q}\left\{\zeta(m_1,t_1)\zeta(m_2,t_2)\cdots\zeta(m_l,t_l) \bigg|\sum_{j=1}^l m_j = w,\quad m_j \geq m,\quad t_j \in T \right\}.$$ \par
	For an irreducible root system $\Phi$, denote $\mathcal{H}(\Phi)$ to be the set of coefficients occurring when writing its highest root in terms of simple roots, it can be easily read off from the following table. \begin{table}[h]
		\begin{tabular}{c|c|c}
			$\Phi$ & $\text{highest root}$     & $\mathcal{H}(\Phi)$                                                        \\ \hline
			$A_n$  & $\alpha_1+\alpha_2+\cdots+\alpha_n$           &  $\{1\}$                                   \\
			$B_n$  & $\alpha_1+2\alpha_2+2\alpha_3+\cdots+2\alpha_n$  & $\{1,2\}$                                \\
			$C_n$  & $2\alpha_1+2\alpha_2+\cdots+2\alpha_{n-1}+\alpha_n$        &  $\{1,2\}$                    \\
			$D_n$  & $\alpha_1+2\alpha_2+\cdots+2\alpha_{n-2}+\alpha_{n-1}+\alpha_n$    & $\{1,2\}$               \\
			$G_2$  & $3\alpha_1+2\alpha_2$     & $\{2,3\}$                                                        \\
			$F_4$  & $2\alpha_1+3\alpha_2+4\alpha_3+2\alpha_4$            & $\{2,3,4\}$                             \\
			$E_6$  & $\alpha_1+2\alpha_2+2\alpha_3+3\alpha_4+2\alpha_5+\alpha_6$      & $\{1,2,3\}$                 \\
			$E_7$  & $2\alpha_1+2\alpha_2+3\alpha_3+4\alpha_4+3\alpha_5+2\alpha_6+\alpha_7$       &  $\{1,2,3,4\}$     \\
			$E_8$  & $2\alpha_1+3\alpha_2+4\alpha_3+6\alpha_4+5\alpha_5+4\alpha_6+3\alpha_7+2\alpha_8$ & $\{2,3,4,5,6\}$
		\end{tabular}
		\caption{\small Highest root of each irreducible root system \cite[Section~12]{humphreys2012introduction}, with respect to some labeling of simple roots.}
		\label{highest_root_table}
	\end{table}

	\begin{theorem}\label{rationality_theorem}
		Let $\Phi$ be an irreducible root system of rank $n$, $r$ its number of positive roots. Let $m$ be a positive integer, we have $$[I_\Phi(s)][(s+m)^{-(r-n)}] \in \sum_{i=1}^n \mathcal{Z}(T)_{m+1,i,mr+n},$$
		where $T$ consists of rational number in $(0,1]$ whose denominator is an element of $\mathcal{H}(\Phi)$.
	\end{theorem}
	
	For a fixed $\Phi$, we let $\delta$ be the number of odd Weyl degrees $d_i$, explicitly,
	$$\delta =   \begin{cases}
		\lfloor \frac{n}{2} \rfloor \quad &\Phi = A_n \text{ for any }n,\\
		1 \quad &\Phi = D_n \text{ with }n\geq 5\text{ and odd},\\
		2 \quad &\Phi = E_6,\\
		0 \quad &\text{otherwise}.\\
	\end{cases}$$
	
	\begin{corollary}\label{coeff_form_corollary}
		Let $T$ as in the statement of Theorem \ref{rationality_theorem}. \\
		(a) (Theorem \ref{special_value_thm}) For even positive integer $m$,
		$$\pi^{mr}\xi_\Phi^{(n)}(-m) \in \sum_{i=1}^n \mathcal{Z}(T)_{m+1,i,mr+n}.$$
		(b) For odd positive integer $m$,
		$$\pi^{mr+n-\delta}\xi_\Phi^{(\delta)}(-m) \in \sum_{i=1}^n \mathcal{Z}(T)_{m+1,i,mr+n}.$$
	\end{corollary}
	\begin{proof}
		For even positive integer $m$, from formula \ref{expr_for_zetaPhi}, we have
		$$\begin{aligned}|W| \times [\xi_\Phi(s)][(s+m)^n] = \left( (2\pi i)^{-m} m!\right)^i [I_\Phi(s)][(s+m)^{-r+n}] \end{aligned}.$$
		For odd positive integer $m$, recall formula \ref{sine_form},
		$$\xi_\Phi(s) \prod_{k=1}^n \frac{\sin(\frac{\pi s d_k}{2})}{\sin(\frac{\pi s}{2})} = \left( \frac{(2\pi)^s}{\Gamma(s)}\right)^r I_\Phi(s),$$
		this implies
		$$\pm \left(\frac{\pi}{2}\right)^{n-\delta} \left(\prod_{d_k \text{ even}} d_k\right) [\xi_\Phi(s)][(s+m)^\delta] = \left( (2\pi)^{-m} m!\right)^r [I_\Phi(s)][(s+m)^{-r+n}].$$
		In both cases, Theorem \ref{rationality_theorem} produces the conclusion. 
	\end{proof}
	
	\begin{remark}
		(1) It is known that $\xi_\Phi(-m)\in \mathbb{Q}$ for all root system $\Phi$ and all $m\in \mathbb{Z}_{\geq 0}$, \cite{bruna2025polynomials}. Therefore part (b) of above corollary gives non-trivial information only for those $\Phi$ with $\delta>0$.  \\
		(2) In the case $\Phi = A_n, B_n, C_n$ or $D_n$, $\mathcal{H}(\Phi)$ equals either $\{1\}$ or $\{1,\frac12\}$, since $\zeta(s,\frac12) = (2^s-1)\zeta(s)$, the numbers in above corollary can be expressed in terms of Riemann's zeta alone.
	\end{remark}
	
	\begin{example}
		Let $\Phi = A_2$, then $\delta = 1$. When $m\in \mathbb{N}$, Corollary \ref{coeff_form_corollary} says
		$$\begin{cases}\pi^{3m}\zeta''_{A_2}(-m) &m\equiv 0 \pmod{2} \\ \pi^{3m+1}\zeta'_{A_2}(-m) &m\equiv 1\pmod{2} \end{cases}$$
		is a $\mathbb{Q}$-linear combination of 
		$$\{\zeta(3m+2), \zeta(i)\zeta(3m+2-i) \bigl | m+1\leq i\leq 2m+1\}.$$
		We will give an explicit form of the coefficients in Remark \ref{A2_alternative_proof} below.
	\end{example}
	
	\begin{example}
		Let $\Phi = B_2$, then $\delta = 0$. When $m\in 2\mathbb{N}$, Corollary \ref{coeff_form_corollary} says
		$$\pi^{4m}\zeta''_{B_2}(-m) \in \text{Span}_\mathbb{Q}\{\zeta(4m+2), \zeta(i)\zeta(4m+2-i) \bigl | m+1\leq i\leq 3m+1\}.$$
	\end{example}
	
	\begin{example}
		Let $\Phi = G_2$, then $\delta = 0$ and $T = \{1,\frac12,\frac13,\frac23\}$, we can remove $\frac23$ because of the relation $\zeta(s,\frac13)+\zeta(s,\frac23) = (3^s-1)\zeta(s)$. When $m\in 2\mathbb{N}$, Corollary \ref{coeff_form_corollary} says
		$$\pi^{6m}\zeta''_{G_2}(-m)$$
		is a $\mathbb{Q}$-linear combination of 
		$$\left\{\zeta(6m+2), \zeta(6m+2,\frac13), \zeta(i)\zeta(6m+2-i),  \zeta(i,\frac13)\zeta(6m+2-i), \zeta(i,\frac13)\zeta(6m+2-i,\frac13)\biggr | m+1\leq i\leq 5m+1\right\}.$$
		The same conclusion holds if one replaces $\zeta(s,1/3)$ with $L(s)$, where $L(s)$ is the non-trivial Dirichlet $L$-function modulo $3$.
	\end{example}
	
	\begin{example}For positive integer $m$, the number $$\begin{cases}\pi^{6m}\zeta^{(3)}_{A_3}(-m) &m\equiv 0 \pmod{2} \\ \pi^{6m+2}\zeta'_{A_3}(-m) &m\equiv 1\pmod{2} \end{cases}$$
		is in $\mathbb{Q}$-span of 
		$$\{\zeta(6m+3), \zeta(k)\zeta(6m+3-k),\zeta(i_1)\zeta(i_2)\zeta(6m+3-i_1-i_2) \big| i_1\geq m+1,i_2\geq m+1,i_1+i_2\leq 5m+2,m+1\leq k\leq 5m+2\}.$$
		When $m$ is even, the numbers $\pi^{9m} \zeta^{(3)}_{B_3}(-m)$ and $\pi^{9m} \zeta^{(3)}_{C_3}(-m)$ are in $\mathbb{Q}$-span of
		$$\{\zeta(9m+3), \zeta(k)\zeta(9m+3-k),\zeta(i_1)\zeta(i_2)\zeta(9m+3-i_1-i_2) \big| i_1\geq m+1,i_2\geq m+1,i_1+i_2\leq 8m+2,m+1\leq k\leq 8m+2\}.$$
	\end{example}
	
	Our proof of Theorem~\ref{rationality_theorem} is qualitative in nature. Although it could, in principle, be made explicit, doing so would require a huge amount of computation, and we therefore do not pursue explicit formulas here.

	\subsection{Proof of Theorem \ref{rationality_theorem}  when $\Phi = A_2$}
	As an illustration, let us prove Theorem \ref{rationality_theorem} for the special case $\Phi = A_2$. We will generalize this procedure in the next two sections for arbitrary $\Phi$.
	
	\begin{lemma}\label{residue_1D_lemma}
		Let $\theta>0$ be a positive real number, let $m \in \mathbb{Z}_{\geq 0}, A\in \mathbb{N}$, $E_s(x)$ be a function analytic in a neighbourhood of $(s,x) \in \{-m\} \times [0,\theta] \subset \mathbb{C}^2$, then $\int_0^\theta x^{As-1} E_s(x) dx$ admits meromorphic continuation to $s\in \mathbb{C}$ and we have
		$$\left[\int_0^\theta x^{As-1} E_s(x) dx\right][(s+m)^{-1}] = \frac{1}{A} [E_{-m}(x)][x^{Am}].$$
	\end{lemma}
	\begin{proof}
		We can write $$\int_0^\theta x^{As-1} E_s(x) dx = \frac{1}{e^{2A \pi i s}-1} \int_{C(\theta)} x^{As-1} E_s(x) dx,$$ with contour $C(\theta)$ in the figure. 	\begin{figure}[h]
			\centering
			\begin{tikzpicture}[decoration={markings,
					mark=at position 0.3cm with {\arrow[line width=1pt]{>}},
					mark=at position 5cm with {\arrow[line width=1pt]{>}},
					mark=at position 7.85cm with {\arrow[line width=1pt]{>}},
					mark=at position 13cm with {\arrow[line width=1pt]{>}}
				}
				]
				\draw[help lines,->] (-1,0) -- (8,0) coordinate (xaxis);
				\draw[help lines,->] (0,-1) -- (0,1) coordinate (yaxis);

				\path[draw,line width=0.8pt,postaction=decorate] (7.5,0) -- (7.5,0.3) -- (0.0,0.3) \centerarcpath(0,0)(90:270:0.3) -- (7.5,-0.3) -- (7.5,0);
				
				\filldraw[black] (7.5,0) circle (1pt) node[anchor=west]{$\theta$};
			\end{tikzpicture}\caption{\small The contour $C(\theta)$.}
		\end{figure}
		The integral on the RHS is analytic around $s=-m$. Therefore
		$$\left[\int_0^\theta x^{As-1} E_s(x) dx\right][(s+m)^{-1}] = \frac{1}{A} \left(\frac{1}{2\pi i}\int_{C(\theta)} x^{-Am-1} E_{-m}(x) dx \right).$$
		The expression in parenthesis, by residue theorem, equals $[E_{-m}(x)][x^{Am}]$.
	\end{proof}
	
	\begin{lemma}\label{zeta_expansion}
		Let $m\geq 1$ and $\Re(a)>0$. As a power series around $x=0$, we have
		$$\zeta(1+m,a+x) = \sum_{k=0}^\infty \binom{m+k}{k} \zeta(1+m+k,a) (-x)^k.$$
	\end{lemma}
	\begin{proof}
		\begin{align*}
			\zeta(1+m,a+x) &= \sum_{n\geq 0} \frac{1}{(n+a+x)^{1+m}} \\ &= \sum_{n\geq 0} \frac{1}{(n+a)^{1+m}} \left(1+\frac{x}{a+n}\right)^{-1-m} \\ &= \sum_{n\geq 0} \sum_{k\geq 0} \frac{1}{(n+a)^{1+m}} \binom{m+k}{k} \left( \frac{-x}{n+a}\right)^k
		\end{align*}
		Exchanging the summation gives the result.
	\end{proof}
	
	Let $m\in \mathbb{N}$, we proceed as follows $$\begin{aligned}I_{A_2}(s) &= \int_0^1 \zeta(1-s,x) \zeta(1-s,1-x)^2 = \int_0^{1/2} \zeta(1-s,x) \zeta(1-s,1-x)^2 + \zeta(1-s,1-x) \zeta(1-s,x)^2 dx \\
		&=  \int_0^{1/2} [\zeta(1-s,1+x)+x^{s-1}] \zeta(1-s,1-x)^2 dx + \int_0^{1/2} \zeta(1-s,1-x) [\zeta(1-s,1+x)+x^{s-1}]^2 dx \\
		&= \int_0^{1/2} x^{s-1} \zeta(1-s,1-x)^2 dx + \int_0^{1/2} 2x^{s-1} \zeta(1-s,1-x)\zeta(1-s,1+x) dx \\ &\quad + \int_0^{1/2} x^{2s-2} \zeta(1-s,1-x) dx + f(s)
	\end{aligned}$$
	where $f(s) = \int_0^{1/2} \zeta(1-s,1+x)\zeta(1-s,1-x)^2 dx + \int_0^{1/2} \zeta(1-s,1-x)\zeta(1-s,1+x)^2 dx$. Because $1\pm x$ are positive on the integration interval $x\in[0,1/2]$, $\zeta(1-s,1\pm x)$ are analytic on $(s,x)\in (\mathbb{C}-\{0\}) \times [0,1/2] \subset \mathbb{C}^2$, therefore $f(s)$ is analytic at $s=-m$ and $[f(s)][(s+m)^{-1}] = 0$. Invoking Lemma \ref{residue_1D_lemma}, 
	$$[I_{A_2}(s)][(s+m)^{-1}] = [\zeta(1+m,1-x)^2][x^m] + [2\zeta(1+m,1-x)\zeta(1+m,1+x)][x^m] + \frac{1}{2}[\zeta(1+m,1-x)][x^{2m+1}],$$
	then using Lemma \ref{zeta_expansion}, which says
	$$\zeta(1+m,1+x) =  \sum_{k=0}^\infty \binom{m+k}{k} \zeta(1+m+k) (-x)^k,$$ yields
	$$[I_{A_2}(s)][(s+m)^{-1}] = \sum_{j=0}^m \binom{m+j}{j}\binom{2m-j}{m-j} \zeta(1+m+j)\zeta(1+2m-j) [1+2(-1)^j] + \frac{1}{2}\binom{3m+1}{2m+1} \zeta(3m+2).$$
	
	Hence we established Theorem \ref{rationality_theorem}  when $\Phi = A_2$. 
	\begin{remark}\label{A2_alternative_proof}
		From $$\xi_{A_2}(s) \frac{(e^{2\pi i s}-1)(e^{3\pi i s}-1)}{(e^{\pi i s}-1)^2} =  \left(\frac{(2\pi i)^s}{\Gamma(s)}\right)^3 I_{A_2}(s),$$ one sees
		$$\begin{cases}\xi_{A_2}''(-m) = \frac{1}{3} ((-2\pi i)^{-m} m!)^3 [I_{A_2}(s)][(s+m)^{-1}] , \quad m\text{ even},\\
			\xi_{A_2}'(-m) = \frac{i}{\pi} ((-2\pi i)^{-m} m!)^3 [I_{A_2}(s)][(s+m)^{-1}], \quad m\text{ odd}, \end{cases}$$
		so we obtain an alternative proof of the formulas of $\xi''_{A_2}(-2\mathbb{N}), \xi'_{A_2}(1-2\mathbb{N})$ given by Onodera \cite{onodera2014functional}. 
	\end{remark}
	
	\section{Proof of Theorem \ref{rationality_theorem} (analytical part)}
	Now we proceed to prove Theorem \ref{rationality_theorem}. First we fix notations that will be used in this section. Let
	\begin{itemize}[leftmargin=*]
		\item $V\subset \mathbb{R}^n$ be an $n$-dimensional simplex, defined over $\mathbb{Q}$ (i.e. all vertices are in $\mathbb{Q}^n$),
		\item $f_1(\boldsymbol{x}),f_2(\boldsymbol{x}),\cdots, g_1(\boldsymbol{x}),g_2(\boldsymbol{x}),\cdots$ be affine linear function on $\mathbb{R}^n$, with rational coefficients, that do not vanish on interior of $V$,
		\item $m$ be a fixed positive integer,
		\item $E_{s}(\boldsymbol{x})$ be a function analytic on a neighbourhood of $(s,\boldsymbol{x})\in \mathbb{C}\times V \subset \mathbb{C}^{n+1}$. 
	\end{itemize}
	
	Lemma \ref{KeyLemma} says
	\begin{equation}\int_V (f_1(\boldsymbol{x})\cdots f_k(\boldsymbol{x}))^{s} E_s(\boldsymbol{x}) d\boldsymbol{x},\qquad \Re(s)>0,\end{equation}
	admits meromorphic extension to $\mathbb{C}$, order of pole at $s=-m$ is at most $n$.
	
	\begin{lemma}\label{V_lemma}
		Let $\mathcal{V}\subset \mathbb{C}$ be a vector space over $\mathbb{Q}$. Assume $E_{-m}(\boldsymbol{x})$ is a polynomial with coefficient in $\mathcal{V}$, then $$\left[\int_V (f_1(\boldsymbol{x})\cdots f_k(\boldsymbol{x}))^{s} E_s(\boldsymbol{x}) d\boldsymbol{x}\right][(s+m)^{-n}] \in \mathcal{V}.$$
	\end{lemma}
	\begin{proof}
		We proceed by induction on $n$. First consider the case $n=1$. Integral (\ref{int_Es}) can be reduced to $$\int_0^1 x^{As} (1-x)^{Bx} f(x)^s E_s(x) dx,$$
		here $A,B \in \mathbb{N}$, $f(x)\in \mathbb{Q}[x]$ non-vanishing on $[0,1]$. Similar to the proof of Lemma \ref{KeyLemma}, one breaks the integration range into two intervals $[0,1/2]$ and $[1/2,1]$. It suffices to prove corresponding statement regarding $$I(s):=\int_0^1 x^{As} f(x)^s E_s(x) dx,$$
		here $A \in \mathbb{N}_{\geq 0}$, $f(x)\in \mathbb{Q}[x]$ non-vanishing on $[0,1]$. If $A=0$, then $I(s)$ is entire, so $[I(s)][(s+m)^{-1}] = 0$. If $A\geq 1$, Lemma \ref{residue_1D_lemma} implies
		$$[I(s)][(s+m)^{-1}] = \frac{1}{A}[f(x)^{-m} E_{-m}(x)][x^{-Am-1}].$$ 
		We assumed $E_{-m}(x)$ has coefficients in $\mathcal{V}$, $f(x)$ is non-vanishing at $x=0$ and has rational coefficients, so above number is also in $\mathcal{V}$. Completing proof for $n=1$ case.\par
		
		We proceed to general $n$, we will closely follow the strategy outlined in Step 3 of Lemma \ref{KeyLemma}, we urge the readers to revisit its proof. \par
		
		Triangulate $V$ using barycentric subdivision, say $V = \cup_i V_i$, the simplexes $V_i$ are defined over $\mathbb{Q}$ since $V$ is.  Write $f(\boldsymbol{x}) := f_1(\boldsymbol{x})\cdots f_k(\boldsymbol{x})$ and let $\mathcal{H}$ consists of zero sets of $f_1,\cdots,f_k$. By hypothesis, they are disjoint from interior of $V$. Write $$\int_V (f(\boldsymbol{x}))^{s} E_s(\boldsymbol{x}) d\boldsymbol{x}= \sum_i \int_{V_i} (f(\boldsymbol{x}))^{s} E_s(\boldsymbol{x}) d\boldsymbol{x},$$
		it suffices to prove, for each $i$, $[\int_{V_i}(f(\boldsymbol{x}))^{s} E_s(\boldsymbol{x}) d\boldsymbol{x}][(s+m)^{-n}] \in \mathcal{V}$. Fix one such $V_i$, we can write $$f(\boldsymbol{x}) = f_{i_1}(\boldsymbol{x}) \cdots f_{i_r}(\boldsymbol{x}) q(\boldsymbol{x}),$$
		where $f_{i_1},\cdots,f_{i_r}$ vanish at the distinguished vertex $p_i$; $q$ is non-zero on $V_i$. Series coefficients of  $E_{i,s}(\boldsymbol{x}) := E_s(\boldsymbol{x}) q(\boldsymbol{x})^s$ when expanded around $\boldsymbol{x} = p_i$ are in $\mathcal{V}$, because we assumed $E_{-m}(\boldsymbol{x})$ satisfies this condition and $q(\boldsymbol{x})$ is a product of rational linear polynomials. Then $$\int_{V_i} (f(\boldsymbol{x}))^{s} E_s(\boldsymbol{x}) d\boldsymbol{x} = \int_{V_i} f_{i_1}(\boldsymbol{x})^s \cdots f_{i_r}(\boldsymbol{x})^s E_{i,s}(\boldsymbol{x}) d\boldsymbol{x}.$$
		
		Translating the point $p_i$ to origin, then applying a suitable linear transformation will map $V_i$ to the standard $n$-simplex 
		$$S_n := \{(x_1,\cdots,x_n) \in \mathbb{R}^n | x_1\geq 0, \cdots, x_n\geq 0, 0\leq x_1+\cdots+x_n\leq 1\}.$$
		
		The Jacobian of this transformation is a rational number since $V_i$ is defined over $\mathbb{Q}$. We reduced the proof into showing that, for \begin{equation}\label{aux_7}I(s) := \int_{S_n} f_1(\boldsymbol{x})^s \cdots f_k(\boldsymbol{x})^s E_s(\boldsymbol{x}) d\boldsymbol{x},\end{equation}
		$[I(s)][(s+m)^{-n}] \in \mathcal{V}$, here $f_1,\cdots,f_k$ are homogeneous linear and $E_s(\boldsymbol{x})$ has series coefficient in $\mathcal{V}$ when expanded around $\boldsymbol{x}=0$. To $I(s)$, we make the same substitution $(x_1,\cdots,x_n) \mapsto (u,y_1,\cdots,y_{n-1})$ as in equation (\ref{int_sub}). Write $\boldsymbol{y} = (y_1,\cdots,y_{n-1})$, we have $f_i(\boldsymbol{x}) = u g_i(\boldsymbol{y})$ with $g_i$ a rational linear polynomial in $\boldsymbol{y}$. Then $$I(s) = \int_{S_{n-1}} g_1(\boldsymbol{y})^s \cdots g_k(\boldsymbol{y})^s  \left(\int_0^1 u^{ks+n-1} E_s(u \boldsymbol{y}, u(1-y_1-\cdots-y_{n-1})) du\right) d\boldsymbol{y}.$$
		
		Let $$\begin{aligned}F_s(\boldsymbol{y}) &:= \frac{e^{2\pi i k s}-1}{2\pi i}\int_0^1 u^{ks + n-1} E_s( u \boldsymbol{y}, u(1-y_1-\cdots-y_{n-1})) du \\ 
			&= \frac{1}{2\pi i} \int_C u^{ks + n-1} E_s( u \boldsymbol{y}, u(1-y_1-\cdots-y_{n-1})) du,\end{aligned}$$
		where $C$ is the contour in Figure \ref{figureC}. By residue theorem, $$F_{-m}(\boldsymbol{y}) = [E_{-m}( u \boldsymbol{y}, u(1-y_1-\cdots-y_{n-1}))][u^{km-n}],$$
		so $F_{-m}(\boldsymbol{y})$ is a polynomial in $y_1,\cdots,y_{n-1}$ with coefficient in $\mathcal{V}$. Finally, 
		$$I(s) = \frac{2\pi i}{e^{2\pi i ks}-1} \int_{S_{n-1}} g_1(\boldsymbol{y})^s \cdots g_k(\boldsymbol{y})^s F_s(\boldsymbol{y}) d\boldsymbol{y},$$
		implies
		$$[I(s)][(s+m)^{-n}] = \frac{1}{k} \left[\int_{S_{n-1}} g_1(\boldsymbol{y})^s \cdots g_k(\boldsymbol{y})^s F_s(\boldsymbol{y}) d\boldsymbol{y}\right][(s+m)^{-n+1}].$$
		Applying the induction hypothesis shows the above number is in $\mathcal{V}$.
	\end{proof}
	
	Let $\mathfrak{V} = \{\mathcal{V}_i\}_{i=0}^\infty$ be an indexed family of $\mathbb{Q}$-vector subspace of $\mathbb{C}$. For an entire function $E(\boldsymbol{x})$ on $\mathbb{C}^n$, we say it is $\mathfrak{V}$-defined at a point $(z_1,\cdots,z_n)$ if 
	$$[E(\boldsymbol{x})][(x_1-z_1)^{k_1}\cdots (x_n-z_n)^{k_n}] \in \mathcal{V}_{\sum k_i}, \qquad \forall (k_1,\cdots,k_n)\in (\mathbb{Z}_{\geq 0})^n.$$
	Observe that if $E(\boldsymbol{x})$ is $\mathfrak{V}$-defined at $\boldsymbol{x}=0$, then 
	$$\frac{1}{2\pi i}\int_C u^{-m-1} E(uy_1,\cdots,uy_{n-1},u(1-y_1-\cdots-y_{n-1})) du,\qquad m\in \mathbb{N},$$
	is a polynomial in $y_1,\cdots,y_{n-1}$ whose coefficients are elements of $\mathcal{V}_m$.
	
	\begin{lemma}\label{V_defined_lemma}
		Assume all $f_1(\boldsymbol{x}),\cdots,f_k(\boldsymbol{x})$ vanish at a vertex $p$ of $V$. If $E_{-m}(\boldsymbol{x})$ is $\mathfrak{V}$-defined at $p$, then $$\left[\int_V f_1(\boldsymbol{x})^{s}\cdots f_k(\boldsymbol{x})^{s} E_s(\boldsymbol{x})d\boldsymbol{x}\right][(s+m)^{-n}] \in \mathcal{V}_{mk-n}.$$
	\end{lemma}
	\begin{proof}
		Again we use induction on $n$. The base case $n=1$ reduces to proving $$\left[\int_0^1 x^{ks} E_s(x) dx \right][(s+m)^{-1}] \in \mathcal{V}_{mk-1},$$
		here $p$ is the origin. By Lemma \ref{residue_1D_lemma}, the LHS equals $\frac{1}{k} [E_{-m}(x)][x^{mk-1}]$, which is in $\mathcal{V}_{mk-1}$ by our assumption on $E_{-m}(x)$. \par 
		For the induction step, we translate $p$ to origin, and apply a suitable rational linear transformation that takes $V$ to the standard simplex $S_n$. It remains to prove the statement for integrals of form $$I(s) := \int_{S_n} f_1(\boldsymbol{x})^s \cdots f_k(\boldsymbol{x})^s E_s(\boldsymbol{x}) d\boldsymbol{x},$$
		where $f_i(\boldsymbol{x})$ are homogeneous, linear, defined over $\mathbb{Q}$; $E_s(\boldsymbol{x})$ analytic on a neighbourhood of $(s,\boldsymbol{x})\in \mathbb{C}\times S_n$ that is $\mathfrak{V}$-defined at origin. We continue the notations used in the proof above, beginning from equation (\ref{aux_7}), 
		$$I(s) = \frac{2\pi i}{e^{2\pi i ks}-1} \int_{S_{n-1}} g_1(\boldsymbol{y})^s \cdots g_k(\boldsymbol{y})^s F_s(\boldsymbol{y}) d\boldsymbol{y},$$
		where $$F_s(\boldsymbol{y}) = \frac{1}{2\pi i} \int_C u^{ks + n-1} E_s( u \boldsymbol{y}, u(1-y_1-\cdots-y_{n-1})) du.$$
		The remark prior to this proposition implies $F_{-m}(\boldsymbol{y})$ is a polynomial with coefficient in $\mathcal{V}_{mk-n}$. Appealing to Lemma \ref{V_lemma} with $\mathcal{V}_{mk-n}$ in place of $\mathcal{V}$ says $(s+m)^{-n+1}$-coefficient of $\int_{S_{n-1}} g_1(\boldsymbol{y})^s \cdots g_k(\boldsymbol{y})^s F_s(\boldsymbol{y}) d\boldsymbol{y}$ lies in $\mathcal{V}$, so
		$$[I(s)][(s+m)^{-n}] = \frac{1}{k} \left[ \int_{S_{n-1}} g_1(\boldsymbol{y})^s \cdots g_k(\boldsymbol{y})^s F_s(\boldsymbol{y})\right]$$
		also lies in $\mathcal{V}$.
	\end{proof}
	
	\begin{lemma}
		Let $g_1(\boldsymbol{x}),\cdots,g_l(\boldsymbol{x})$ be positive on $V$; $f_1,\cdots,f_k$ be non-negative on $V$, and they all vanish at a vertex $p$ of $V$. Let $T = \{g_1(p),\cdots,g_l(p)\}$, then $$\left[\int_V f_1(\boldsymbol{x})^{s-1}\cdots f_k(\boldsymbol{x})^{s-1} \zeta(1-s,g_1(\boldsymbol{x}))\cdots \zeta(1-s,g_l(\boldsymbol{x}))d\boldsymbol{x} \right] \left[(s+m)^{-n}\right] \in \mathcal{Z}(T)_{m+1,l,(m+1)(k+l)-n}.$$ 
	\end{lemma}
	\begin{proof}
		Let $E_s(\boldsymbol{x}) = (1+s)^l\zeta(-s,g_1(\boldsymbol{x}))\cdots \zeta(-s,g_l(\boldsymbol{x}))d\boldsymbol{x}$, it is analytic on a neighbourhood of $(s,\boldsymbol{x}) \in \mathbb{C}\times V$ because $g_i(\boldsymbol{x})$ are positive on $V$. Denote the family of vector space $\mathfrak{V} = \{\mathcal{V}_i\}_{i=0}^\infty$ by $\mathcal{V}_i = \mathcal{Z}(T)_{m+1,l,l(m+1)+i}$. By Lemma \ref{zeta_expansion}, 
		$$\zeta(1+m,a+x) = \sum_{k=0}^\infty \binom{m+k}{k} \zeta(1+m+k,a) (-x)^k, \qquad a>0,$$
		we see $E_{-m-1}(\boldsymbol{x})$ is $\mathfrak{V}$-defined at point $\boldsymbol{x}=p$. The lemma above says $(s+m+1)^{-n}$-coefficient of
		$$\left[\int_V f_1(\boldsymbol{x})^{s}\cdots f_k(\boldsymbol{x})^{s} E_s(\boldsymbol{x})d\boldsymbol{x}\right] \left[(s+m+1)^{-n} \right] \in \mathcal{V}_{(m+1)k-n},$$ which is equivalent to our statement.
	\end{proof}
	
	Our main result of this section is the following.
	\begin{proposition}\label{zetacoeffform}
		Let $g_1(\boldsymbol{x}),\cdots,g_N(\boldsymbol{x})$ be non-negative on $V$, and suppose they do not simultaneously vanish at any vertex of $V$. Let $$T = \left\{g_i(p) | 1\leq i \leq N, p \text{ is a vertex of }V\right\}\cup\{1\} - \{0\}.$$ Then, 
		$$\left[\int_V \zeta(1-s,g_1(\boldsymbol{x})) \cdots \zeta(1-s,g_N(\boldsymbol{x})) d\boldsymbol{x}\right] \left[(s+m)^{-n}\right] \in \sum_{i=1}^{N-n} \mathcal{Z}(T)_{m+1,i,(m+1)N-n}.$$
	\end{proposition}
	\begin{proof}
		We triangulate $V$ using barycentric subdivision, say $V = \bigcup_i V_i$, each $V_i$ is $\mathbb{Q}$-defined. It suffices to prove, for each $i$, $\int_{V_i} \zeta(1-s,g_N(\boldsymbol{x})) \cdots \zeta(1-s,g_N(\boldsymbol{x}))$ satisfies the statement. Recall $V_i$ comes with a distinguished vertex $p_i$ of $V$, such that $\forall j, g_j$ is either non-vanishing on $V_i$ or $g_j(p_i) = 0$. Fix an $i$, after some re-numbering, we can assume $g_1,\cdots,g_l$ vanish at $p_i$ ; $g_{l+1},\cdots,g_N$ non-zero on $V_i$. Using $\zeta(1-s,x) = \zeta(1-s,x+1) + x^{s-1}$, write above integral as
		\begin{equation}\label{aux_4}\int_{V_i} [g_{1}(\boldsymbol{x})^{s-1} + \zeta(1-s,1+g_{1}(\boldsymbol{x}))]\cdots [g_{l}(\boldsymbol{x})^{s-1} + \zeta(1-s,1+g_{l}(\boldsymbol{x}))] \zeta(1-s,g_{l+1}(\boldsymbol{x})) \cdots \zeta(1-s,g_N(\boldsymbol{x}))  d\boldsymbol{x}.\end{equation}
		Expand this into $2^{l}$ terms, each term is of form
		\begin{equation}\label{aux_0}\int_{V_i} g_{a_1}(\boldsymbol{x})^{s-1} \cdots g_{a_k}(\boldsymbol{x})^{s-1} \zeta(1-s,1+g_{a_{k+1}}(\boldsymbol{x})) \cdots \zeta(1-s,1+g_{a_{l}}(\boldsymbol{x})) \zeta(1-s,g_{l+1}(\boldsymbol{x})) \cdots \zeta(1-s,g_N(\boldsymbol{x}))  d\boldsymbol{x},\end{equation}
		here $\{a_1,\cdots,a_l\}$ is a permutation of $\{1,\cdots,l\}$, $k$ ranges from $0$ to $l$ inclusive. This last integral satisfies requirement of previous lemma: $g_{a_1},\cdots,g_{a_k}$ vanishes at the vertex $p_i$ of $V_i$; $1+g_{a_{k+1}},\cdots,1+g_{a_l}, g_{l+1},\cdots,g_N$ are positive on $V_i$, so $(s+m)^{-n}$ coefficient of expression (\ref{aux_0}) is an element of $\mathcal{Z}(T)_{m+1,N-k,(m+1)N-n}$. Observe that $k \leq n-1$ does not contribute: by Corollary \ref{orderofpolelessthank}, integral (\ref{aux_0}) have poles of order $\leq k\leq n-1$. Therefore $(s+m)^{-n}$-coefficient of expression (\ref{aux_4}) is an element of 
		$$\sum_{k=n}^l \mathcal{Z}(T)_{m+1,N-k,(m+1)N-n},$$
		now recall the definition of $l$ (which depends on the sub-simplex $V_i$): it is the number of linear forms among $\{g_1(\boldsymbol{x}),\cdots,g_l(\boldsymbol{x})\}$ that vanish at $p_i$. Since we assumed $g_1(\boldsymbol{x}),\cdots,g_N(\boldsymbol{x})$ do not simultaneously vanish at any vertex of $V$, $l \leq N-1$. \par
		
		Summing over all simplexes $V_i$ that triangulate $V$, we have
		$$\left[\int_V \zeta(1-s,g_1(\boldsymbol{x})) \cdots \zeta(1-s,g_N(\boldsymbol{x})) d\boldsymbol{x}\right] \left[(s+m)^{-n}\right] \in \sum_{k=n}^{N-1} \mathcal{Z}(T)_{m+1,N-k,(m+1)N-n},$$
		as desired.
	\end{proof}
	
	Proposition \ref{zetacoeffform} actually almost implies Theorem \ref{rationality_theorem}, we illustrate this for the root system $B_2$.
	
	\begin{example}\label{B2_example} Recall the expression for $I_{B_2}(s)$ from equation (\ref{int_rep_I}), 
		$$I_{B_2}(s) = \int_{[0,1]^2} \zeta(1-s,x_1)\zeta(1-s,x_2) \zeta(1-s,\{-x_1-2x_2\})\zeta(1-s,\{-x_1-x_2\}) dx_1dx_2$$
		To remove fractional part, we rewrite the integral into four parts (as shown in the figure):
		\begin{multline*}I_{B_2}(s) = \int_{W_1} \zeta(1-s,x_1)\zeta(1-s,x_2) \zeta(1-s,1-x_1-2x_2)\zeta(1-s,1-x_1-x_2) dx_1dx_2 \\
			+ \int_{W_2} \zeta(1-s,x_1)\zeta(1-s,x_2) \zeta(1-s,2-x_1-2x_2)\zeta(1-s,1-x_1-x_2) dx_1dx_2 \\
			+ \int_{W_3} \zeta(1-s,x_1)\zeta(1-s,x_2) \zeta(1-s,2-x_1-2x_2)\zeta(1-s,2-x_1-x_2) dx_1dx_2 \\
			+ \int_{W_4} \zeta(1-s,x_1)\zeta(1-s,x_2) \zeta(1-s,3-x_1-2x_2)\zeta(1-s,2-x_1-x_2) dx_1dx_2.
		\end{multline*}
		\vspace*{-7mm}
		\begin{figure}[h]
			\centering
			\scalebox{1}{
				\begin{tikzpicture}
					\tikzmath{\x1 = -4; \y1 =-0.9;
						\x2 = 1.6; \y2 = -1.2;
						\x3 = 0; \y3 = 4; 
						\x4 = (\x1 + \x2)/2; \y4 = (\y1 + \y2)/2;
						\x5= (\x1 + \x3)/2; \y5= (\y1 + \y3)/2;
						\x6 = (\x2 + \x3)/2; \y6 = (\y2 + \y3)/2;}
					\draw (0,0)--(0,6);
					\draw (0,6)--(6,6);
					\draw (6,6)--(6,0);
					\draw (6,0)--(0,0);
					\draw (0,3)--(6,0);
					\draw (6,0)--(0,6);
					\draw (0,6)--(6,3);
					\filldraw[black] (0,0) circle (1pt) node[anchor=east]{$(0,0)$};
					\filldraw[black] (6,6) circle (1pt) node[anchor=west]{$(1,1)$};
					\filldraw[black] (6,3) circle (1pt) node[anchor=west]{$(1,\frac{1}{2})$};
					\filldraw[black] (6,0) circle (1pt) node[anchor=west]{$(1,0)$};
					\filldraw[black] (0,3) circle (1pt) node[anchor=east]{$(0,\frac{1}{2})$};
					\filldraw[black] (0,6) circle (1pt) node[anchor=east]{$(0,1)$};
					\filldraw[black] (2.5,0.5) node[anchor=south]{$W_1$};
					\filldraw[black] (1.5,3) node[anchor=south]{$W_2$};
					\filldraw[black] (4.5,3) node[anchor=north]{$W_3$};
					\filldraw[black] (4.3,4.5) node[anchor=south]{$W_4$};
			\end{tikzpicture}}
			\caption{\small Triangulation of square $[0,1]^2$ for $I_{B_2}(s)$.}
		\end{figure}
		
		We apply Proposition \ref{zetacoeffform} to each of the four integrals. For example, for the integral over $W_1$, we take $(g_1,g_2,g_3,g_4) = (x_1,x_2,1-x_1-2x_2,1-x_1-x_2)$, the $T$ is simply the values of $g_i$ attained at vertices of $W_1$, one easily calculates $T = \{1,\frac{1}{2}\}$. Repeat this for $W_2, W_3, W_4$, we arrive at
		$$[I_{B_2}(s)][(s+m)^{-2}] \in \sum_{i=1}^2 \mathcal{Z}(\{1,\frac{1}{2}\})_{m+1,i,4(m+1)-2},$$
		therefore we proved Theorem \ref{rationality_theorem} in the case $\Phi = B_2$.
	\end{example}
	
	To prove Theorem \ref{rationality_theorem} in its full generality, we need to find a triangulation of $[0,1]^{r-n}$ that removes the fractional part in $I_\Phi(s)$ with controlled "denominator". This is the goal of the next section.
	
	\section{Proof of Theorem \ref{rationality_theorem} (combinatorial part)}
	Let $A\in \text{Mat}_{n\times n}(\mathbb{Z})$ be invertible, define its level, denoted by $\text{level}(A)$, as the smallest $N\in \mathbb{N}$ such that $NA^{-1} \in \text{Mat}_{n\times n}(\mathbb{Z})$. For $1\leq j\leq n$, we let $\varepsilon_j\in \mathbb{R}^n$ to be the vector with all entries $0$ except the $j$-th entry, which is $1$. Some easy observations:
	\begin{enumerate}
		\item Level is invariant under transpose, that is $\text{level}(A) = \text{level}(A^T)$.
		\item If $U_1,U_2 \in \text{GL}_n(\mathbb{Z})$, then $\text{level}(A) = \text{level}(U_1 AU_2)$.
		\item Permuting columns (resp. rows) of $A$ does not change its level.
		\item If $A$ has a column that is of form $\varepsilon_j$, let $X$ be the $(n-1)\times (n-1)$ matrix obtained by removing this column and the $j$-row of $A$, then $\text{level}(X) = \text{level}(A)$.
	\end{enumerate}
	Here $(4)$ follows because of the identity $$\begin{pmatrix}1 & v \\ 0 & X \end{pmatrix}^{-1} = \begin{pmatrix}1 & -vX^{-1} \\ 0 & X^{-1} \end{pmatrix}.$$
	
	The following lemma is also evident.
	\begin{lemma}\label{lemma_level}
		Let $A$ be an invertible $n\times n$ matrix with integer entries, let $y\in \mathbb{Z}^n$. The unique solution $x\in \mathbb{Q}^n$ of $Ax = y$ satisfies $\text{level}(A) x \in \mathbb{Z}^n$.
	\end{lemma}
	
	For an $n\times m$ integral matrix $B$ with $m>n$, define $$\mathcal{D}(B) := \{\text{level}(A) | A \text{ ranges over all invertible } n\times n \text{ sub-matrices of }A\} \subset \mathbb{N},$$
	note that there are $\binom{m}{n}$ ways of forming $n\times n$ sub-matrices from $B$.

	\begin{lemma}\label{transpose_D}
		Let $C$ be an $M\times N$ integral matrix, $I_N$ be the identity matrix of size $N$, then $\mathcal{D}$ of the following two block matrices coincide:
		$$\begin{pmatrix}I_M & C \end{pmatrix}, \qquad \begin{pmatrix}I_N & C^T \end{pmatrix}.$$
	\end{lemma}
	\begin{proof}
		We denote $B_1 := \begin{pmatrix}I_M & C \end{pmatrix}$ and $B_2 :=  \begin{pmatrix}I_N & C^T \end{pmatrix}$. Let $A$ be an $M\times M$ sub-matrix of $B_1$, write $$A = \begin{pmatrix}\varepsilon_{i_1} & \cdots & \varepsilon_{i_k} & v_1 & \cdots & v_{M-k} \end{pmatrix},$$
		here $v_1,\cdots,v_{M-k}$ are columns from $C$. Then by applying repeatedly observation $(4)$ above, we see $\text{level}(A)$ equals to the level of the $(M-k)\times (M-k)$ matrix obtained from $\begin{pmatrix} v_1 & \cdots & v_{M-k} \end{pmatrix}$ after removing the rows with indices $i_1,\cdots,i_k$. Therefore
		$$\mathcal{D}(B_1) = \{\text{level}(X) | X \text{ ranges over all invertible sub-matrices of any size of }C\}.$$
		By replacing $C$ with $C^T$, we have
		$$\mathcal{D}(B_2) = \{\text{level}(X) | X \text{ ranges over all invertible sub-matrices of any size of }C^T\}.$$
		Because transpose of sub-matrices are sub-matrices of transpose, and levels are invariant under transpose, we have $\mathcal{D}(B_1) = \mathcal{D}(B_2)$.
	\end{proof}
	
	The next lemma is crucial in what follows, it says $\mathcal{D}(B)$ are essentially "denominators" of the triangulation coming from columns of $B$.
	
	\begin{lemma}\label{tri_B_lemma}
		Let $B$ be an $n\times m$ integral matrix with $m>n$, assume $B$ contains the identity matrix $I_n$ as a sub-matrix. Let $l_i: \mathbb{R}^n \to \mathbb{R}$ be the linear function coming from\footnote{more precisely, if $i$-th column of $B$ is $(a_1,\cdots,a_n)$, then $l_i(x_1,\cdots,x_n) = a_1 x_1+\cdots + a_nx_n$.} $i$-th column of $B$. Then there exists a triangulation of the hypercube $[0,1]^n = \cup_j V_j$, such that for each $j$,
		\begin{enumerate}
			\item coordinates of each vertex of $V_j$ are rational numbers whose denominator divides some element of $\mathcal{D}(B)$.
			\item for each $i$, there exists $N_{i,j}\in \mathbb{Z}$ such that $N_{i,j}\leq l_i(\boldsymbol{x}) \leq N_{i,j}+1$ for all $\boldsymbol{x} \in V_j$.
		\end{enumerate}
	\end{lemma}
	\begin{proof}
		For any $n$ linearly independent columns of $B$, corresponds to linear functions $l_{i_1}(\boldsymbol{x}),\cdots,l_{i_n}(\boldsymbol{x})$. Let $\mathcal{P}_{\{i_1,\cdots,i_n\}} \subset \mathbb{R}^n$ be the collection of $\boldsymbol{x} \in \mathbb{Q}^n$ such that 
		\begin{equation}\label{aux_2}l_{i_1}(\boldsymbol{x})\in \mathbb{Z}, \cdots, l_{i_n}(\boldsymbol{x})\in \mathbb{Z}.\end{equation}
		Since $l_{i_1},\cdots,l_{i_n}$ are linearly independent, for each choices of integers, above system of linear equation is uniquely solvable, so $\mathcal{P}_{\{i_1,\cdots,i_n\}}\subset \mathbb{R}^n$ is discreet, its intersection with compact set $[0,1]^n$ is thus finite. For any $n$-element subset $I$ of $\{1,2,\cdots,m\}$, we define $\mathcal{P}_I$ as above if the corresponding columns are linearly independent, and $\mathcal{P}_I = \varnothing$ otherwise. Define 
		$$\mathcal{P} := \bigcup_{\substack{ I \subset \{1,\cdots,m\} \\ \# I = n}} \mathcal{P}_I \cap [0,1]^n \subset \mathbb{Q}^n \cap [0,1]^n.$$
		This set is finite since each $\mathcal{P}_I \cap [0,1]^n$ is finite. Each coordinate of $\boldsymbol{x}\in \mathcal{P}$ is a rational number whose denominator divides some element of $\mathcal{D}(A)$ because it solves equation (\ref{aux_2}), see Lemma \ref{lemma_level}. Also, since $B$ contains the identity matrix, the $2^n$ vertices of the hypercube $[0,1]^n$ are in $\mathcal{P}$. \par
		Define a set $\mathfrak{T}$, whose elements are of form $\{W_1,\cdots,W_k\}$ such that
		\begin{itemize}
			\item $W_1,\cdots,W_k$ are $n$-simplexes,
			\item $\{W_1,\cdots,W_k\}$ triangulate $[0,1]^n$,
			\item vertices of each $W_i$ are in $\mathcal{P}$.
		\end{itemize}
		Since $\mathcal{P}$ contains $2^n$ vertices of the hypercube, $\mathfrak{T}$ is non-empty. Since $\mathcal{P}$ is finite, only finitely many simplexes can be formed with vertices in $\mathcal{P}$, so $\mathfrak{T}$ is a finite set. Pick an element $\{V_1,\cdots,V_M\}$ of $\mathfrak{T}$ with maximal $M$. We claim that $\{V_1,\cdots,V_M\}$ satisfies the requirement of the lemma. \par
		Point (1) is satisfied by construction. We only need to check point (2). Suppose not, then for some $i, j$ and integer $N$, the hyperplane $H := \{\boldsymbol{x}: l_i(\boldsymbol{x}) = N\}$ intersects the interior\footnote{here $\text{span}(v_0,\cdots,v_n)$ means convex span of points $v_0,\cdots,v_n$, which is the $n$-simplex with these point as vertices.} of $V_j := \text{span}(v_0,\cdots,v_n)$. Their intersection $H\cap V_j$ is an $(n-1)$-simplex, spanned by $n$ points, say $c_1,\cdots,c_n$. 
		\begin{figure}[h]
			\centering
			\includegraphics[width=0.4\textwidth]{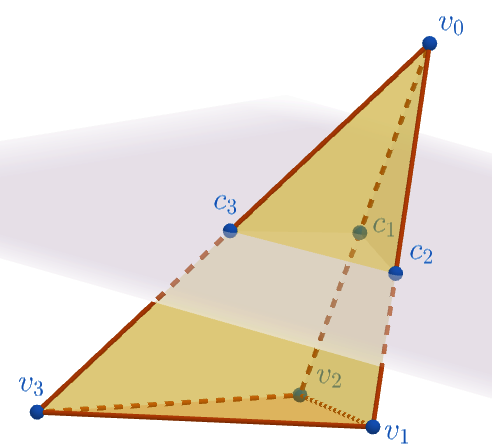}
			\caption{\small Illustrative case when $n=3$.}
		\end{figure} 
		Each $c_i\in \mathcal{P}$ because all faces of $V_j$ are of the form $l_k(\boldsymbol{x})\in \mathbb{Z}$. It follows that we can further triangulate (see figure for illustration) $$V_j = W_1 \cup W_2 \cup \cdots W_N, \qquad N\geq 2,$$ into smaller pieces, all vertices of $W_i$ are also in $\mathcal{P}$. We now produced another element of $\mathfrak{T}$, namely 
		$$\{V_1,\cdots,V_{j-1},W_1,\cdots,W_N,V_{j+1},\cdots,V_M\},$$
		whose cardinality is greater than $M$, this contradicts our choice of $M$. Hence (2) must hold.
	\end{proof}
	
	For a root system $\Phi$ of rank $n$. Fix a set of positive roots $\Phi^+ = \{\alpha_1,\cdots,\alpha_r\}$, let $\{\lambda_1,\cdots,\lambda_n\}$ be the associated fundamental weights. Consider the $n\times r$ integral matrix $M_\Phi := (\lambda_i,\alpha_j^\vee)$, we abuse notation and define $\mathcal{D}(\Phi)$ to be $$\mathcal{D}(\Phi):= \mathcal{D}(M_\Phi) \subset \mathbb{N}.$$
	
	\begin{proposition}\label{root_system_T_prop}
		For an irreducible root system $\Phi$, $m \in \mathbb{N}$, we have $$[I_\Phi(s)][(s+m)^{-(r-n)}]\in \sum_{k=1}^n \mathcal{Z}(T)_{m+1,k,mr+n},$$
		where $T$ consists of rational numbers in $(0,1]$ whose denominator divides an element of $\mathcal{D}(\Phi)$.
	\end{proposition}
	\begin{proof}
		Consider variables $(x_\alpha)$ indexed by $\alpha \in \Phi^+-\Delta$, there are $r-n$ of them, consider the $(r-n)\times r$ matrix formed by taking coefficient of linear forms
		$$\left\{x_\alpha | \alpha \in \Phi^+-\Delta \right\} \quad \text{ and } \quad \left\{\sum_{\alpha \in \Phi^+ - \Delta} (\lambda_i, \alpha^\vee) x_\alpha \biggr| 1\leq r\leq n\right\}.$$
		Denote this matrix by $B$. Note that the $n\times r$ matrix $M_\Phi$ contains the identity matrix $I_n$ as sub-matrix, write it in the block form $M_\Phi = \begin{pmatrix}I_n & C \end{pmatrix}$ for some integral matrix $C$, then $B = \begin{pmatrix} I_{r-n} & C^T \end{pmatrix}$, so Lemma \ref{transpose_D} implies $\mathcal{D}(B) = \mathcal{D}(M_\Phi) = \mathcal{D}(\Phi)$. \par
		Apply Lemma \ref{tri_B_lemma} to the $(r-n)\times r$ matrix $B$, we obtain a triangulation of $[0,1]^{r-n} = \cup_j V_j$, and
		$$\begin{aligned}I_\Phi(s) &=  \int_{[0,1]^{r-n}} \prod_{\alpha \in \Phi^+ - \Delta} \zeta(1-s,x_\alpha) \prod_{i=1}^n \zeta\left(1-s,\left\{ -\sum_{\alpha \in \Phi^+ - \Delta} (\lambda_i, \alpha^\vee) x_\alpha \right\}\right) dx_\alpha \\
			&=  \sum_j \int_{V_j} \prod_{\alpha \in \Phi^+ - \Delta} \zeta(1-s,x_\alpha) \prod_{i=1}^n \zeta\left(1-s,\left\{ -\sum_{\alpha \in \Phi^+ - \Delta} (\lambda_i, \alpha^\vee) x_\alpha \right\}\right) dx_\alpha  \\ &:= \sum_j I_{V_j}(s). \end{aligned}$$
		Point (2) of this lemma says $V_j$ resolves the fractional part of the above integrand, and point (1) says coordinate of each vertices of $V_j$ is divisible by an element of $\mathcal{D}(B)$. Therefore values of the linear functions $$x_\alpha, \qquad \left\{-\sum_{\alpha \in \Phi^+ - \Delta} (\lambda_i, \alpha^\vee) x_\alpha \right\},$$
		on $V_j$'s vertices lie in $T$. Invoking Proposition \ref{zetacoeffform} gives $$[I_{V_j}(s)][(s+m)^{-(r-n)}] \in \sum_{k=1}^n \mathcal{Z}(T)_{m+1,k,mr+n}.$$
		Summing this over $j$ completes the proof.
	\end{proof}
	
	To prove Theorem \ref{rationality_theorem}, it remains to find $\mathcal{D}(\Phi)$ for each irreducible root system $\Phi$. 
	
	\begin{example}\label{exception_D_ex}
		For $\Phi = G_2$, we have $$M_\Phi = \begin{pmatrix}
			1& 0& 1& 1& 2& 1 \\ 0& 1& 3& 1& 3& 2
		\end{pmatrix}.$$
		After computing levels of all $\binom{6}{2} = 15$ two-by-two submatrices, we find $\mathcal{D}(G_2) = \{1,2,3\}$. This proves Theorem \ref{rationality_theorem} for $\Phi = G_2$. \par
		
		For $\Phi = B_3$, we have $$M_\Phi= \begin{pmatrix}
			1 & 0 & 0 & 1 & 0 & 0 & 1 & 2 & 1 \\
			0 & 1 & 0 & 1 & 1 & 2 & 1 & 2 & 2 \\
			0 & 0 & 1 & 0 & 1 & 1 & 1 & 1 & 1 \\
		\end{pmatrix}.$$
		After exhaustively computing levels of all $\binom{9}{3} = 84$ three-by-three submatrices, we find $\mathcal{D}(B_3) = \{1,2\}$. This proves Theorem \ref{rationality_theorem} for $\Phi = B_3$. 
	\end{example}
	
	\subsection{Determining $\mathcal{D}(\Phi)$}
	For a (not necessarily irreducible) root system $\Phi$ in the Euclidean space $V$. Let $L(\Phi) := \text{Span}_\mathbb{Z} \Delta$ be its root lattice. We fix a base $\Delta = \{\alpha_1,\cdots,\alpha_n\}$ consists of simple roots. 
	
	Let $S \subset \Phi$ such that $S$ spans $V$, the lattice generated by $S$ has full rank, so its quotient with the root lattice has finite index. We define
	\begin{equation}\label{def_E}\mathcal{E}(\Phi) := \left\{\text{Exponent of the finite abelian group }\frac{L(\Phi)} {\text{Span}_\mathbb{Z} S} \biggr| S\subset \Phi \text{ such that }\text{Span}_\mathbb{R}S=V\right\} \subset \mathbb{N}.\end{equation}
	
	\begin{proposition}\label{DE_prop}
		We have $\mathcal{D}(\Phi) = \mathcal{E}(\Phi^\vee)$, where $\Phi^\vee$ is the dual root system of $\Phi$.
	\end{proposition}
	\begin{proof}
		In the definition of $\mathcal{E}(\Phi)$, we can assume $\#  S=n, S\in \Phi^+$, say $S = \{\beta_1,\cdots,\beta_n\} \subset \Phi^+$. Such $S$ is in bijective correspondence with invertible $n\times n$ matrices of $M_\Phi$. Let $A = (\lambda_i, \beta_j^\vee)$ be the corresponding $n\times n$ sub-matrix of $M_\Phi$. \par
		Fix an $1\leq i \leq n$, let $(c_{i1},\cdots,c_{in})$ be the $i$-th column of $A^{-1}$, then $A(c_{i1},\cdots,c_{in})^T = \varepsilon_i$. Writing $y = \sum_{k=1}^n c_{ik} \beta_k^\vee$, we see $$\begin{cases}
			(\lambda_j, y) &= 0 \qquad  i\neq j \\
			(\lambda_j,y) &= 1 \qquad  i=j
		\end{cases}.$$
		This forces $y = \alpha_i^\vee$ because fundamental weights are dual to simple co-roots. The level of $A$ is the smallest positive integer $N$ such that $N c_{ij} \in \mathbb{Z}, 1\leq i,j\in n.$  On the other hand, $\alpha_j^\vee = \sum_{k=1}^n c_{jk} \beta_k^\vee$ says $N$ is the smallest positive integer such that $\{N\alpha_1^\vee,\cdots,N\alpha_n^\vee\}\subset \text{Span}_\mathbb{Z}(S^\vee)$, that is, the exponent of the abelian group
		$$\frac{L(\Phi^\vee)} {\text{Span}_\mathbb{Z} S^\vee},$$
		as desired.
	\end{proof}
	
	\begin{proposition}\label{EH_prop}
		For an irreducible root system $\Phi$, $\mathcal{E}(\Phi) = \mathcal{H}(\Phi) \cup \{1\}$.
	\end{proposition}
	
	\begin{remark}\label{springer_remark}
		In \cite[Chapter~1,~Section 4]{Springer1970}, a weaker result is recorded, namely, primes $p$ for which $\frac{L(\Phi)} {\text{Span}_\mathbb{Z} S}$ can have non-trivial $p$-torsion are exactly the primes in $\mathcal{H}(\Phi)$. 
	\end{remark}
	
	On inclusion in Proposition \ref{EH_prop} is obvious: $\mathcal{H}(\Phi) \cup \{1\} \subset \mathcal{E}(\Phi)$. Indeed, $1\in \mathcal{E}(\Phi)$ by letting $S= \Delta$ in the definition (\ref{def_E}). If $q$ appears as a coefficient of the highest root $\alpha_0$, say $$\alpha_0 = q\alpha_i + \text{(linear combination of }\alpha_j \text{ with } i\neq j),$$
	letting $S = \{\alpha_0\} \cup \{\alpha_j | j\neq i\}$, then it is evident that
	$$\frac{L(\Phi)} {\text{Span}_\mathbb{Z} S} \cong \frac{\mathbb{Z}}{q\mathbb{Z}} \implies q\in \mathcal{E}(\Phi). $$
	
	Next we proceed to establish the non-trivial inclusion: $\mathcal{E}(\Phi) \subset \mathcal{H}(\Phi) \cup \{1\}$. For a subset $I$ of $\Delta$, we let $V_I$ be the Euclidean subspace of $V$ spanned by $I$, $\Phi_I := \Phi \cap V_I$. Note that $\Phi_I$ is a root system with base $I$ \cite[Chapter~1,~Section~10]{humphreys1992reflection}.
	
	\begin{lemma}
		Let $\beta_1,\cdots,\beta_k$ linear independent roots of $\Phi$, then there exists some $\sigma\in W$ and $I\subset \Delta$, such that $\# I =k$ and  $\{\sigma(\beta_1),\cdots,\sigma(\beta_k)\} \subset \Phi_I.$
	\end{lemma}
	\begin{proof}
		Let $V'$ be the $\mathbb{R}$-span of $\beta_1,\cdots,\beta_k$, then $\Phi' := \Phi \cap V'$ is a sub-root system of $\Phi$, it has rank $k$. Let $\{\alpha'_1,\cdots,\alpha'_k\}$ be a base of $\Phi'$. Base of a sub-root system can always be extended to a base of the ambient root system \cite[Chapter~6,~Proposition~24]{bourbaki2002lie}, so there exists $\{\alpha'_{k+1},\cdots,\alpha'_n\} \subset \Phi$ such that
		$$\{\alpha'_1,\cdots,\alpha'_k,\alpha'_{k+1},\cdots,\alpha'_n\}$$
		is a base of $\Phi$. Since $W$ acts transitively on bases \cite[Chapter~10]{humphreys2012introduction}, there exists $\sigma\in W$ such that $$\{\sigma(\alpha'_1),\cdots,\sigma(\alpha'_k),\sigma(\alpha'_{k+1}),\cdots,\sigma(\alpha'_n)\} = \{\alpha_1,\cdots,\alpha_n\},$$
		say $\sigma(\alpha'_1) = \alpha_{i_1},\cdots, \sigma(\alpha'_k) = \alpha_{i_k}$. Let $I = \{\alpha_{i_1}, \cdots, \alpha_{i_k}\}$. Then for $1\leq j\leq k$, $$\sigma(\beta_j) \in \sigma(\Phi') \subset \sigma(\text{span}_\mathbb{R}(\alpha'_1,\cdots,\alpha'_k)) = V_I.$$
		Therefore $\sigma(\beta_j) \in \Phi \cap V_I = \Phi_I$, as desired.
	\end{proof}
	
	The following lemma says $\mathcal{E}(\Phi)$ can be computed inductively, removing one simple root at a time.
	\begin{lemma}
		Let $\Phi$ be an irreducible root system, $\alpha_0$ be its highest root. Then $$\mathcal{E}(\Phi) \subset \bigcup_{i=1}^n \mathcal{E}_i(\Phi),$$ where \begin{equation}\label{aux_9} \mathcal{E}_i(\Phi) := \left\{\text{LCM}(q,e) | e\in \mathcal{E}(\Phi_{\Delta-\{\alpha_i\}}), 1\leq q\leq \text{(coefficient of }\alpha_i \text{ in }\alpha_0 \text{)}\right\}.\end{equation}
	\end{lemma}
	\begin{proof}
		Let $S = \{\beta_1,\cdots,\beta_n\} \subset \Phi$ be linear independent, applying above lemma to $\{\beta_1,\cdots,\beta_{n-1}\}$, we obtain a Weyl group element $\sigma$, a permutation $\pi: \{1,\cdots,n\}\to \{1,\cdots,n\}$ and some $(n-1)\times (n-1)$ integral matrix such that
		$$\begin{pmatrix} \sigma(\beta_1) \\ \vdots \\ \sigma(\beta_n)\end{pmatrix} = \begin{pmatrix} A & 0 \\ \ast & q\end{pmatrix} \begin{pmatrix} \alpha_{\pi(1)} \\ \vdots \\ \alpha_{\pi(n)}\end{pmatrix},\quad q\in \mathbb{N}.$$
		Because $S$ spans $V$, $q\neq 0$, replacing $\beta_n$ by its negative if necessary, we can assume $q>0$. Now, with $I = \Delta - \{\alpha_{\pi(n)}\}$, we have
		\begin{equation}\label{aux_8}\frac{L(\Phi)} {\text{Span}_\mathbb{Z}S} \cong \frac{\mathbb{Z}}{q\mathbb{Z}} \bigoplus \frac{\mathbb{Z}^{n-1}} {A\mathbb{Z}^{n-1}} \cong \frac{\mathbb{Z}}{q\mathbb{Z}} \bigoplus \frac{L(\Phi_I)} {\text{Span}_\mathbb{Z}S'},\end{equation}
		where $S' = \{\sigma(\beta_1),\cdots,\sigma(\beta_{n-1})\} \subset \Phi_I$. The number $q$ is the $\alpha_{\pi(n)}$-coefficient of $\sigma(\beta_n)\in \Phi$, hence is less or equal to the $\alpha_{\pi(n)}$-coefficient in $\alpha_0$. Exponent of (\ref{aux_8}) is the LCM of $q$ and the exponent of $\frac{L(\Phi_I)} {\text{Span}_\mathbb{Z}S'}$, the latter is an element of $\mathcal{E}(\Phi_{\Delta-\{\alpha_{\pi(n)}\}})$.
	\end{proof}
	
	If $\Phi$ is reducible, say $\Phi = \Phi_1\times \Phi_2$, then it is obvious from definition (\ref{def_E}) that
	$$\mathcal{E}(\Phi) = \{\text{LCM}(a,b) | a\in \mathcal{E}(\Phi_1), b\in \mathcal{E}(\Phi_2)\}.$$
	
	For irreducible $\Phi$, we proceed to show $\mathcal{E}(\Phi) \subset \mathcal{H}(\Phi) \cup \{1\}$ by induction. For each $1\leq i\leq n$, we consider $\Phi_I$ with $I = \Delta-\{\alpha_i\}$, this amounts to deleting one vertex in its Dynkin diagram, and by induction, $\mathcal{E}(\Phi_{\Delta-\{\alpha_i\}})$ can be considered as known, so is $\mathcal{E}_i(\Phi)$ in equation (\ref{aux_9}).
	\begin{figure}[h]
		\centering 
		\begin{tikzpicture}
			\dynkin[text style/.style={scale=1}, labels={1,2,3,4},
			label directions={,,,,,},
			edge length=.75cm,
			scale=1.8, mark=o] F4
		\end{tikzpicture} \qquad \qquad 
		\begin{tikzpicture}
			\dynkin[text style/.style={scale=1}, labels={1,2,3,4,5,6},
			label directions={,,,,,},
			edge length=.75cm,
			scale=1.8, mark=o] E6
		\end{tikzpicture} \\
		\begin{tikzpicture}
			\dynkin[text style/.style={scale=1}, labels={1,2,3,4,5,6,7},
			label directions={,,,,,},
			edge length=.75cm,
			scale=1.8, mark=o] E7
		\end{tikzpicture}
		\begin{tikzpicture}
			\dynkin[text style/.style={scale=1}, labels={1,2,3,4,5,6,7,8},
			label directions={,,,,,},
			edge length=.75cm,
			scale=1.8, mark=o] E8
		\end{tikzpicture}
		\caption{\small Dynkin diagrams of $F_4, E_6, E_7, E_8$, with labeling as in Table \ref{highest_root_table}.}
		\label{exceptional_Dynkin}
	\end{figure}
	
	The base case to start is $\Phi = A_1$, in which $\mathcal{E}(A_1) = \{1\}$. Recall the explicit expression of highest roots given in Table \ref{highest_root_table}. For each $1\leq i\leq n$, we write $q$ and $e$ as in equation (\ref{aux_9}) and $I = \Delta - \{\alpha_i\}$.
	
	\begin{itemize}[leftmargin=*]
		\item $\Phi = A_n$. For any $i$, the coefficient of $\alpha_i$ in the highest root is $1$, so $q=1$. Also, $\Phi_I$ is isomorphic to some $A_j\times A_k$ with $j,k<n$, so $e=1$ by induction hypothesis. Therefore $\mathcal{E}_i(\Phi) \subset \{1\} \implies \mathcal{E}(\Phi) = \{1\}$.
		\item $\Phi = X_n, X\in  \{B,C,D\}$. For any $i$, we have $q\in \{1,2\}$ and $\Phi_I$ is isomorphic to a product of $X_j (j<n)$ and some copies of $A_k$, so $e\in \{1,2\}$ by induction hypothesis. Therefore $\mathcal{E}_i(\Phi) \subset \{1,2\} \implies \mathcal{E}(\Phi) \subset \{1,2\}$.
		\item $\Phi = F_4$. The highest root is $2\alpha_1+3\alpha_2+4\alpha_3+2\alpha_4$. \begin{itemize} \item When $i=1$, $q\in \{1,2\}$. Since $\Phi_I = C_3$, $e\in \{1,2\}$ by the $C_3$ case above. Therefore $\mathcal{E}_1(\Phi) \subset \{1,2\}$. 
			\item When $i=2$, $\Phi_I = A_2\times A_1$, and $q\leq 3$, so $\mathcal{E}_2(\Phi) \subset \{1,2,3\}$. 
			\item When $i=3$, $\Phi_I = A_2\times A_1$ and $q\leq 4$, $\mathcal{E}_3(\Phi) \subset \{1,2,3,4\}$. 
			\item When $i=4$, $q\in \{1,2\}$. Since $\Phi_I = B_3$, $e\in \{1,2\}$ by the $B_3$ case above. Therefore $\mathcal{E}_4(\Phi) \subset \{1,2\}$. 
		\end{itemize}
		\item $\Phi = E_6$. The highest root is $\alpha_1+2\alpha_2+2\alpha_3+3\alpha_4+2\alpha_5+\alpha_6$. \begin{itemize} \item When $i=1$ or $6$, $\Phi_I = D_5$, $q\in \{1,2\}$ and $e\in \{1,2\}$ by the $D_5$ case. So $\mathcal{E}_1(\Phi) \subset \{1,2\}$. 
			\item When $i=2$, $\Phi_I = A_5$, and $q\leq 2$, $\mathcal{E}_2(\Phi) \subset \{1,2\}$
			\item When $i=3$ or $5$, $\Phi_I = A_4\times A_1$, and $q\leq 2$, $\mathcal{E}_3(\Phi) \subset \{1,2\}$
			\item When $i=4$, $\Phi_I = A_2\times A_2\times A_1$, and $q\leq 3$, $\mathcal{E}_4(\Phi) \subset \{1,2,3\}$.
		\end{itemize}
		\item $\Phi = E_7$. The highest root is $2\alpha_1+2\alpha_2+3\alpha_3+4\alpha_4+3\alpha_5+2\alpha_6+\alpha_7$. \begin{itemize} \item When $i=1$, $\Phi_I = D_6$, $q\in \{1,2\}$ and $e\in\{1,2\}$ by the $D_6$ case, $\mathcal{E}_1(\Phi) \subset \{1,2\}$. 
			\item When $i\in\{2,3,4,5\}$, $\Phi_I$ decomposes as products of $A_j$, so $e=1$. Also $q\in\{1,2,3,4\}$, so $\mathcal{E}_i(\Phi) \subset \{1,2,3,4\}$.
			\item When $i=6$, $\Phi_I = D_5\times A_1$, $q\leq 2$ and $e\in \{1,2\}$, so $\mathcal{E}_6(\Phi) \subset \{1,2\}$.
			\item When $i=7$, $\Phi_I = E_6$, $q = 1$ and $e\in \{1,2,3\}$ by the $E_6$ case. Hence $\mathcal{E}_7(\Phi) \subset \{1,2,3\}$.
		\end{itemize}
		\item $\Phi = E_8$. The highest root is $2\alpha_1+3\alpha_2+4\alpha_3+6\alpha_4+5\alpha_5+4\alpha_6+3\alpha_7+2\alpha_8$. One checks exactly as before, we only mention some interesting cases. \begin{itemize} \item When $i=4$, $\Phi_I = A_4\times A_2\times A_1$, $q\leq 6$, so $\mathcal{E}_4(\Phi) \subset \{1,2,3,4,5,6\}$.
			\item When $i=5$, $\Phi_I = A_4\times A_3$, $q\leq 5$, so $\mathcal{E}_5(\Phi) \subset \{1,2,3,4,5\}$.
			\item When $i=7$, $\Phi_I = E_6\times A_1$, $q\leq 3$ and $e\in \{1,2,3\}$ so $\mathcal{E}_7(\Phi) \subset \{1,2,3,4,5,6\}$.
		\end{itemize}
	\end{itemize}
	
	We checked $\mathcal{E}(\Phi)\subset \mathcal{H}(\Phi)\cup \{1\}$ for all irreducible $\Phi$, and Proposition \ref{EH_prop} is thus established.

	Finally we can complete the proof of our second main result.
	\begin{proof}[Proof of Theorem \ref{rationality_theorem}]
		From Proposition \ref{root_system_T_prop}, we have $$[I_\Phi(s)][(s+m)^{-(r-n)}]\in \sum_{k=1}^n \mathcal{Z}(T)_{m+1,k,mr+n},$$
		where $T$ consists of rational numbers in $(0,1]$ with denominator divisible by an element of $\mathcal{D}(\Phi)$. Combining Lemma \ref{DE_prop} and Proposition \ref{EH_prop}, we have $\mathcal{D}(\Phi) = \mathcal{H}(\Phi^\vee) \cup \{1\}= \mathcal{H}(\Phi) \cup \{1\}$. 
	\end{proof}
	
	\section*{Declarations}
	The author declares no conflict of interest. No data were generated or analyzed in this study.

	\bibliographystyle{plain} 
	\bibliography{../ref.bib} 
	
\end{document}